\documentclass[11pt,a4paper]{article}%
\usepackage{amsmath}
\usepackage{graphicx}
\usepackage{epsfig}%
\usepackage{amsfonts}%
\usepackage{amssymb}
\usepackage{makeidx}         
\usepackage{color}
\usepackage{bm}
\usepackage{soul}
\usepackage{appendix}
\usepackage{subfig}
\RequirePackage{mathrsfs}
\usepackage{yfonts}

    \newcommand    {\be}     { \beta  }

    \renewcommand  {\phi}    { \varphi }
    
    \newcommand    {\bq}   { \begin{equation} }
    \newcommand    {\eq}   { \end{equation} }
    \newcommand    {\ba}   { \begin{array} }
    \newcommand    {\ea}   { \end{array} }
\newcommand{\refeq}[1]{Eq.~(\ref{#1})}

\def\diag{\mathop{\rm diag}\nolimits}

\def\max{\mbox{max}}

\def\u{\mnm{u}}
\def\U{\mnm{U}}

\def\thetab{\text{\boldmath $\theta$}}

\def\u{\mnm{u}}
\def\U{\mnm{U}}

\newcommand{\BIGOP}[1]{\mathop{\mathchoice%
{\raise-0.22em\hbox{\huge $#1$}} {\raise-0.05em\hbox{\Large $#1$}}
{\hbox{\large $#1$}}{#1}}}

\newcommand{\BIGboxplus}{\mathop{\mathchoice%
{\raise-0.35em\hbox{\huge $\boxplus$}}%
{\raise-0.15em\hbox{\Large $\boxplus$}}{\hbox{\large
$\boxplus$}}{\boxplus}}}

\usepackage{tikz} \usetikzlibrary{positioning,shapes}
\usepackage{pgfplots} \usetikzlibrary{plotmarks}
\usetikzlibrary{calc,backgrounds}

\def\cov{\mathop{\rm cov}\nolimits}

\def\trace{\mbox{trace}}

\def\sinc{\mbox{sinc}}
\DeclareMathAlphabet{\mathitbf}{OML}{cmm}{b}{it}
\def\thetab{\text{\boldmath $\theta$}}

\def\z{\mathitbf{z}}

\def\u{\mathitbf{u}}

\def\F{\mathitbf{F}}

\def\U{\mathitbf{U}}
\def\H{\mathcal{H}}

\def\jmath{j}

\def\diag{\mathop{\rm diag}\nolimits}
\def\cov{\mathop{\rm cov}\nolimits}

\def\ee{\mathbf{\check{e}}}

\def\trace{\mathop{\rm trace}\nolimits}
\def\det{\mathop{\rm det}\nolimits}

\def\LL{\mathcal{L}}

\def\bse{\begin{eqnarray*}}
\def\ese{\end{eqnarray*}}
\def\be{\begin{eqnarray}}
\def\ee{\end{eqnarray}}
\def\bq{\begin{equation}}

\newcommand{\bB}{\mathbf{B}}

\newcommand{\bC}{\mathbf{C}}

\newcommand{\bL}{\mathbf{L}}

\newcommand{\bW}{\mathbf{W}}

\newcommand{\bx}{\mathbf{x}}
\newcommand{\bz}{\mathbf{z}}

\newcommand{\bu}{\mathbf{u}}
\newcommand{\bv}{\mathbf{v}}

\newcommand{\bw}{\mathbf{w}}
\newcommand{\by}{\mathbf{y}}

\newtheorem{theorem}{Theorem}[section]
\newtheorem{ex}{Example}[section]

\newtheorem{lemma}[theorem]{Lemma}

\newtheorem{remark}[theorem]{Remark}

\makeatletter\@addtoreset{equation}{section}\makeatother
\makeatletter\@addtoreset{figure}{section}\makeatother
\makeatletter\@addtoreset{table}{section}\makeatother
\newcommand{\refeqp}[1]{(Eq.~\ref{#1})}

\begin{document}

\title{Tucker tensor analysis of \\
Mat\'ern functions in spatial statistics}


\author{Alexander Litvinenko
  \and
  David Keyes
  \and
  Venera Khoromskaia
  \and
  Boris N. Khoromskij
  \and
  Hermann G. Matthies
}

\newcommand{\Addresses}{{
  \bigskip
  \footnotesize

  Alexander Litvinenko (Corresponding author), \textsc{King Abdullah University of Science and Technology, 
Thuwal 23955-6900, Kingdom of Saudi Arabia}, 
  \textit{E-mail:} \texttt{alexander.litvinenko@kaust.edu.sa}

  \medskip

  David Keyes, \textsc{King Abdullah University of Science and Technology, 
Thuwal 23955-6900, Kingdom of Saudi Arabia},
  \textit{E-mail:} \texttt{david.keyes@kaust.edu.sa}

  \medskip

  Venera Khoromskaia, \textsc{Max-Planck Institute for Mathematics in the Sciences, D-04103, Leipzig, Germany; Max-Planck Institute for Dynamics of Complex Technical Systems, D-39106, Magdeburg, Germany},
  \textit{E-mail:} \texttt{vekh@mis.mpg.de}
  
  \medskip
  
  Boris Khoromskij, \textsc{Max-Planck Institute for Mathematics in the Sciences, D-04103, Leipzig, Germany},
  \textit{E-mail:} \texttt{bokh@mis.mpg.de}

  \medskip
  
  Hermann G. Matthies, \textsc{Technische Universitaet Braunschweig, D-38106, Braunschweig, Germany},
  \textit{E-mail:} \texttt{wire@tu-braunschweig.de}

}}

\maketitle

\Addresses
          
\maketitle

\begin{abstract}
In this work, we describe advanced numerical tools for working with multivariate functions and for the 
analysis of large data sets. These tools will drastically reduce the required computing time and the 
storage cost, and, therefore, will allow us to consider much larger data sets or finer meshes. 
Covariance matrices are crucial in spatio-temporal statistical tasks, but are often very expensive to 
compute and store, especially in 3D. Therefore, we approximate covariance functions by cheap surrogates 
in a low-rank tensor format. We apply the Tucker and canonical tensor decompositions to a family of 
Mat\'ern- and Slater-type functions with varying parameters and demonstrate numerically 
that their approximations exhibit exponentially fast convergence. 
We prove the exponential convergence of the Tucker and canonical approximations in tensor
rank parameters.
Several statistical operations are performed in this low-rank tensor format, including evaluating the 
conditional covariance matrix, spatially averaged
estimation variance, computing a quadratic form, determinant, trace, loglikelihood, inverse, 
and Cholesky decomposition of a large covariance matrix. 
Low-rank tensor approximations reduce the computing and storage costs essentially. 
For example, the storage cost 
is reduced from an exponential $\mathcal{O}(n^d)$ to a linear scaling $\mathcal{O}(drn)$, 
where $d$ is the spatial dimension, $n$ is the number of mesh points in one direction, 
and $r$ is the tensor rank. 
Prerequisites for applicability of the proposed techniques are the assumptions that the data, locations, 
and measurements lie on a tensor (axes-parallel) grid and that the covariance 
function depends on a distance, $\| x-y\|$. 
\end{abstract}
\textit{\textbf{Corresponding author:}} Dr. Alexander Litvinenko, alexander.litvinenko$@$kaust.edu.sa

\section{Introduction}
Nowadays it is very common to work with large spatial data sets 
\cite{sun2016statistically, Furrer2011, Stein2004, LitvGenton17, Stein2012, nychka2015multiresolution}, 
for instance, with 
satellite data, collected over a very large area (e.g., the data collected by 
the National Center for Atmospheric Research, USA\footnote{https://www.earthsystemgrid.org/}). 
This data can also come from a computer simulator 
code as a solution of a certain multiparametric equation (e.g., Weather research and Forecasting 
model\footnote{https://www.mmm.ucar.edu/weather-research-and-forecasting-model}), 
it could be also sensor data from multiple sources. 
Typical operations in spatial statistics, such as evaluating the 
spatially averaged
estimation variance, computing quadratic forms of the conditional covariance 
matrix, or computing maximum of likelihood function \cite{Stein2004} require high computing 
power and time. Our motivation for using low-rank tensor techniques is that operations on advanced matrices, 
such as hierarchical, low-rank and sparse matrices,
are limited by their high computational costs, especially in 3D and for a large number of observations. 

A tensor can be simply defined as a high-order matrix,
where multi-indices are used instead of indices (see Section~\ref{sec:low-rank} and Eq.~\ref{Tensor_def} 
for a rigorous definition). One way to obtain a tensor from a vector or matrix is to reshape it. For example, 
we assume $\bv \in \mathbb{R}^{10^6}$ is a vector. We reshape it and obtain a matrix of size $10^3 \times 10^3$, 
or a tensor of order $3$ of size $10^2\times 10^2\times 10^2$ or a tensor of order 6 of 
size $10\times ... \times 10$ (6 times). 
Each element of such a six-dimensional hypercube is described by the multi-index $\alpha=(\alpha_1,...,\alpha_6)$. 
The obtained tensors contain not only rows and columns, but also \textit{slices} and \textit{fibers} 
\cite{Kolda:01,Kolda:07,DMV-SIAM2:00}. 
These slices and fibers can be analyzed for linear dependences, super symmetry, or sparsity 
and may result in a strong data compression. 
Another difference between tensors and matrices is that a matrix (obtained, for instance, 
after the discretization 
of a kernel $c(\bx,\by)=c(\vert \bx-\by\vert)$) separates a point $\bx=(x_1,...,x_d)\in \mathbb{R}^{d}$ from 
a point $\by=(y_1,...,y_d)\in \mathbb{R}^{d}$, whereas the corresponding tensor (depending on tensor format) 
separates $x_1-y_1$ from $x_2-y_2$... and $x_d-y_d$. This implies that tensors may have not just one rank like a 
matrix, but many. Therefore, we speak about a \textit{tensor rank}, but not a matrix rank. 
In this work, we consider 
two very common tensor formats: canonical (denoted as CP) and Tucker (see Section~\ref{sec:low-rank}). 

Low-rank tensor mathods can be gainfully combined with other data-compression techniques in low
dimensions. 
For example, a 3D function can be approximated as the sum of the tensor products of one-dimensional function. 
Then the usual matrix techniques can be applied to those 1D functions.

To be more concrete, we consider a relatively wide class of Mat\'{e}rn covariance functions. We demonstrate 
how to approximate Mat\'{e}rn covariance matrices in a low-rank tensor format, then how to perform typical 
kriging and spatial statistics operations in this tensor format. Mat\'ern covariance matrices typically
depend on three to five unknown hyper-parameters, 
such as smoothness, three covariance lengths (in a 3D anisotropic case), and variance. 
We study the dependences of tensor ranks and approximation errors on these parameters. Splitting the spatial 
variables via low-rank techniques reduces the
computing cost for a matrix-vector product from $\mathcal{O}(N^2)$ to $\mathcal{O}(dr^2n^2)$ FLOPs, 
where $d$ is the spatial dimension, $r$ is the tensor rank, and 
$n$ is the number of mesh points along the longest edge 
of the computational domain. For simplicity, we assume that $N=n^d$ (e.g., $d=4$ for a time-space problem in 3D).
Other motivating factors for applying low-rank tensor techniques include the following:
\begin{enumerate}
\item The storage cost is reduced from $\mathcal{O}(n^d)$ to $\mathcal{O}(drn)$ or, depending 
on the tensor format, to $\mathcal{O}(drn+r^d)$, where $d>1$. 
\item The low-rank tensor technique allows us to compute not only the matrix-vector product, 
but also the inverse $\bC^{-1}$, square root $\bC^{1/2}$, 
matrix exponent $\exp(\bC)$, $\trace(\bC)$, $\det(\bC)$, and a likelihood function.
\item The low-rank tensor approximation is relatively new, but already a well-studied technique 
with free software libraries available.
\item The approximation accuracy is fully controlled by the tensor rank. The full rank gives an exact 
representation;
\item Low-rank techniques are either faster than a Fourier transform ($\mathcal{O}(drn)$ 
vs. $\mathcal{O}(n^d \log n^d)$) or 
can be efficiently combined with it \cite{nowak2013kriging, DoKhSa-qtt_fft-2011};
\end{enumerate}
General limitations of the tensor technique are that 1) it could be time consuming to compute a low-rank tensor 
decomposition; 2) it requires axes-parallel mesh; 3) some theoretical estimations exist for functions depending 
on $\vert x-y \vert$ (although more general functions have a low-rank representation in practice).

During the last few years, there has been great interest in numerical methods for 
representing and approximating large covariance matrices
\cite{LitvGenton17, Rasmussen05,saibaba2012application,
nowak2013kriging,ambikasaran2013large,BallaniKressner, LitvHLIBpro}.
Low-rank tensors were previously applied to accelerated kriging and spatial design by orders of 
magnitude \cite{nowak2013kriging}. 
The covariance matrix under consideration was assumed to be circulant, and the 
first column had a low-rank decomposition. Therefore, $d$-dimensional Fourier 
was applied to and drastically reduce the storage and the computing cost. 
 
The maximum likelihood estimator was computed for parameter-fitting given Gaussian observations 
with a Mat\'ern covariance matrix \cite{Minden16}. The presented framework for unstructured 
observations in two spatial dimensions allowed for an evaluation of the log-likelihood 
and its gradient with computational complexity $\mathcal{O}(n^{3/2})$. 
 {The $\H$-matrix techniques \cite{Part1,Part2,HackHMEng} provide the efficient data sparse 
approximation for the differential and integral operators in $\mathbb{R}^d$, $d=1,2,3$}. 
$\H$-matrices are very robust for approximating the covariance matrix 
\cite{khoromskij2009application, saibaba2012application, BallaniKressner, harbrecht2015efficient},
\cite{ambikasaran2013large, BoermGarcke2007}, 
its inverse \cite{ambikasaran2013large}, and its Cholesky decomposition 
\cite{khoromskij2009application, LitvGenton17, LitvHLIBpro}, but can also be expensive, 
especially for large $n$ in 3D.  { Namely, the complexity in 3D will be $C\, k^{d-1} N\log N$, 
where $N=n^d$, $d=3$, $k\ll n$ is the rank and $C $ is a large constant which scales exponentially 
in dimension $d$, see \cite{Part2}. 
Thus, the ${\cal H}$-matrix techniques scale  exponentially in dimension size.}
Therefore, more efficient methods for fast and efficient matrix linear algebra operations are still needed.
 
The key idea is to compute a low-rank decomposition not of the covariance function (it could be hard), 
but of its 
analytically known spectral density (which could be a much easier object) and then apply the inverse 
Fourier to the obtained low-rank components.
The Fourier transformation of the Mat\'{e}rn 
covariance function is known analytically as the Hilbert tensor. 
This Hilbert tensor can be decomposed 
numerically in a low-rank tensor format. Both the Fourier Transformation and its 
inverse have the canonical (CP) tensor rank-1. Therefore, the inverse Fourier does not change the 
tensor rank of the argument. By applying the inverse Fourier to the low-rank tensor, 
we obtain a low-rank approximation of the initial covariance 
matrix, which can be 
further used in the Kalman Filter update, Karhunen-Lo{\`e}ve Expansion, Bayesian Update, 
and kriging.

The structure of the paper is as follows.
In Section \ref{sec:Motivation}, we list typical tasks from statistics that motivate us to use 
low-rank tensor techniques and
define the Mat\'{e}rn covariance functions and their Fourier transformations.
Section \ref{sec:low-rank} is devoted to low-rank tensor decomposition.
Sub-sections~\ref{ssect: Sinc_C2T}, \ref{ssec:Laplace_transf}, and \ref{ssec:Cover_tensor} 
contain the main theoretical 
contribution of this work. We present low-rank tensor techniques and separate radial basis 
functions using the Laplace transform and
the $\sinc$ quadrature, give estimations of the approximation error, convergence rate, and 
the tensor rank, and we also 
prove the existence of a low-rank approximation of a Mat\'{e}rn function.  
Section~\ref{sec:8tasksLow} contains another important contribution of this work, namely, 
the solutions to typical 
statistical tasks in the low-rank tensor format.

\section{Motivation}
\label{sec:Motivation}
\subsection{Problem settings in spatial statistics}
Below, we formulate \textbf{five tasks}. These computational tasks are very common and important in statistics. 
Fast and efficient solution of these tasks will help to solve many real-world problems, such as the weather 
prediction, moisture modeling, and optimal design in geostatistics.\\  
\tikzstyle{decision}=[diamond, draw, fill=blue!30]
\tikzstyle{line} = [draw, -latex']
\tikzstyle{line2} = [draw, -stealth, thick]
\tikzstyle{elli} = [draw, ellipse, fill=red!=50, minimum height=15mm]
\tikzstyle{block} = [draw, rectangle, fill=blue!30, text width=15em, text centered, 
minimum height=16mm, node distance =10em]
\tikzstyle{shortblock} = [draw, rectangle, fill=blue!30, text width=10em, text centered, 
minimum height=16mm, node distance =10em]
\tikzstyle{vshortblock} = [draw, rectangle, fill=blue!30, text width=7em, text centered, 
minimum height=12mm, node distance =7em]
\tikzstyle{evshortblock} = [draw, rectangle, fill=blue!30, text width=4em, text centered, 
minimum height=12mm, node distance =7em]
\tikzstyle{arrow} = [draw, rectangle, fill=white!50, text width=3em, text centered, 
minimum height=7mm, node distance =7em]
%
%
%
\textbf{Task 1. Approximate a Mat\'{e}rn covariance function in a low-rank tensor format.}\\
The covariance function $c(\bx,\by)$, $\bx=(x_1,...,x_d)$,  $\by=(y_1,...,y_d)$, is discretized 
on a tensor grid with $N$ mesh points, $N=n^d$, $d\geq 1$ and $\varepsilon > 0$. The task is 
to find the following decomposition into one-dimensional functions, i.e.,
$\Vert c(\bx, \by) -  \sum_{i=1}^r \prod_{\mu=1}^d c_{i \mu}(x_{\mu}, y_{\mu}) \Vert \leq \varepsilon $, 
for some given $\varepsilon>0$. 
Alternatively, we may look for factors $\bC_{i\mu}$ such that
$\Vert \bC- \sum_{i=1}^r \bigotimes_{\mu=1}^d \bC_{i \mu} \Vert \leq \varepsilon $. 
Here, the matrices $ \bC_{i \mu}$ correspond to the one-dimensional covariance functions 
$c_{i \mu}(x_{\mu}, y_{\mu})$ in the direction $\mu$. \\
\textbf{Task 2. Computing of square root of $\bC$.}
The square root $\bC^{1/2}$  of the covariance matrix $\bC$ is needed in order to generate 
random fields and processes. It is also used in the Kalman filter update.\\
\textbf{Task 3. Kriging.} Spatial statistics and kriging \cite{Kitanidis_book} are used to model 
the distribution of ore 
grade, forecast of rainfall intensities, moisture, temperatures, or contaminant. 
The missing values are interpolated
from the known measurements by kriging \cite{Matheron_1971,Journel_Huijbregts_1978_book}.
When considering space-time kriging on fine meshes
\cite{Wesson_Pegram_2004,finke2004mapping,haylock2008european,  de2011toward},
kriging may easily exceed the computational power of modern computers.
Estimating the variance of kriging and 
geostatistical optimal design problems are especially numerically intensive 
\cite{Mueller_2007_Spatial_Data,Nowak_2010MG_Measures_Uncertainty,Spoeck_Pilz_2010SERRA_OD_convex}.

The kriging can be defined as follows.
Let $\hat{\mathbf{s}}$ be the $N\times1$ vector of values to be estimated
with zero expectation and covariance matrix $\mathbf{C}_{ss}$.
Let $\mathbf{y}$ be the $m\times1$, $m\ll N$, vector of
measurements. The corresponding cross- and auto-covariance matrices
are denoted by $\mathbf{C}_{sy}$ and $\mathbf{C}_{yy}$, and
sized $N\times m$ and $m\times m$, respectively. If the measurements are subject to
error, an error covariance matrix $\mathbf{R}$ is included in $\mathbf{C}_{yy}$.
Using this notation, the kriging estimate $\hat{\mathbf{s}}$ is given by
$\hat{\mathbf{s}}=\mathbf{C}_{sy}\mathbf{C}_{yy}^{-1}\mathbf{y}$.\\
%
%
\textbf{Task 4. Geostatistical design.} The goal of geostatistical optimal design is to 
optimize the sampling patterns
from which the data values in $\mathbf{y}$ will be obtained.
The objective function that will be minimized is typically a scalar 
measure of either the conditional covariance matrix or
the estimation variance \refeqp{eq:estvar}. The two most common measures for 
geostatistical optimal design are $\phi_A$ and $\phi_C$:
%
\begin{equation}
\phi_A  =  N^{-1}\trace\left[\mathbf{C}_{ss|y}\right],\quad \text{ and}\quad 
\phi_C = \mathbf{z}^T(\mathbf{C}_{ss|y})\mathbf{z}   ,
\label{eq:Ccriterion}
\end{equation}
where
$\mathbf{C}_{ss|y}:= \mathbf{C}_{ss} - \mathbf{C}_{sy} \mathbf{C}_{yy}^{-1}\mathbf{C}_{ys}$
\cite{Mueller_2007_Spatial_Data,Nowak_2010MG_Measures_Uncertainty}.\\
%
%
%
\textbf{Task 5. Computing the joint Gaussian log-likelihood function.} We assume 
that $\z \in \mathbb{R}^{N}$ is an available vector of measurements, and $\thetab$ is an unknown vector of 
the parameters of a covariance matrix $\bC$. The task is to compute the maximum likelihood estimation (MLE), 
where the log-likelihood function is as follows
\begin{equation}
\label{eq:like}
\LL(\thetab)=-\frac{N}{2}\log(2\pi) - \frac{1}{2}\log \det\{\bC(\thetab)\}-
\frac{1}{2}(\bz^T\cdot \bC(\thetab)^{-1}\bz ). 
\end{equation}
The difficulty here is that each iteration step of a maximization procedure requires the 
solution of a linear system $\bL \bv=\bz$, the Cholesky 
decomposition, and the determinant.

In Section~\ref{sec:8tasksLow} we give detailed solutions. We give strict definition of 
tensors later in Section~\ref{sec:low-rank}. 

\subsection{Mat\'{e}rn covariance and its Fourier transform}
\label{sec:Matern}
A low-rank approximation of the covariance function is a key component of the tasks formulated above.
Among of the many covariance models available, the Mat\'{e}rn family \cite{Matern1986a,Handcock1993a} 
is widely used in 
spatial statistics, 
geostatistics  \cite{chiles2009geostatistics}, machine learning  \cite{BoermGarcke2007}, image analysis, 
weather forecast, moisture modeling, and as the correlation for temperature fields \cite{North2011a}. 
The work \cite{Handcock1993a} introduced the Mat\'{e}rn form of spatial correlations into statistics 
as a flexible parametric 
class with one parameter determining the smoothness of the underlying spatial random field.

The main idea of this low-rank approximation is shown on the Diagram 1 and explained in 
details in Section~\ref{sec:idea}. 
Diagram 1 demonstrates two possible ways to find a low-rank tensor approximation of the 
Mat\'{e}rn covariance function. 
The first way (marked with `?') is not so trivial and the second via the Fast Fourier 
Transform (FFT), low-rank and the 
inverse FFT (IFFT) is more trivial. We use here the fact that the FT of the Mat\'{e}rn 
covariance is analytically known 
and has a known low-rank approximation. 
The IFFT can be computed numerically and does not change the tensor ranks.

\begin{tikzpicture}
\node[shortblock](process0){$\approx \sum_{i=1}^R \bigotimes_{\mu=1}^d \bC_{i \mu} $};
\node[shortblock, below of=process0](process1){$\approx \sum_{i=1}^R \bigotimes_{\nu=1}^d \tilde{\u}_{i \nu}$};
\node[block, left of=process0, xshift=-12em, yshift=-0em]
(process2){$C_{\nu,\ell}(r)=\frac{2^{1-\nu}}{\Gamma(\nu)}\left(\frac{\sqrt{2\nu}r}{\ell}
\right)^\nu K_\nu\left(\frac{\sqrt{2\nu}r}{\ell}\right)$};
\node[block, left of=process1, xshift=-12em, yshift=-0em](process3){$f_{\alpha,\nu}(\rho)
:=\frac{\Gamma(\nu+d/2)\alpha^{2\nu}}{\pi^{d/2}\Gamma(\nu)}\frac{1}{(\alpha^2+\rho^2)^{\nu+d/2}}$};
\path[line2](process2)--node[yshift=1em]{?}(process0);
\path[line2](process1)--node[xshift=1.4em]{IFFT}(process0);
\path[line2](process3)--node[yshift=1em]{Low-rank }(process1);
\path[line2](process3)--node[yshift=-1em]{approximation}(process1);
\path[line2](process2)--node[xshift=1.4em]{FFT}(process3);
\node [below=1cm, xshift=12em, align=flush center,text width=14cm] at (process3)
{Diagram 1: Two possible ways to find a low rank tensor approximation of 
the Mat\'{e}rn covariance matrix $C_{\nu,\ell}(r)$. 
};
\end{tikzpicture}

The Mat\'{e}rn covariance function is defined as
\begin{equation}
\label{eq:Matern}
C_{\nu,\ell}(r)=\frac{2^{1-\nu}}{\Gamma(\nu)}\left(\frac{\sqrt{2\nu}r}{\ell}
\right)^\nu K_\nu\left(\frac{\sqrt{2\nu}r}{\ell}\right),
\end{equation}
where 
distance $r:=\Vert x-y\Vert $, $x,y$ two points in $\mathbb{R}^d$,
$\nu>0$ defines the smoothness of the random field. The larger $\nu$, the smoother random field.
The parameter $\ell>0$ is a spatial range parameter that measures how quickly the correlation of the random 
field decays with distance, with larger $\ell$ corresponding to a slower decay (keeping $\nu$ fixed).  
$K_\nu$ denotes the modified Bessel function of order $\nu$. 
When $\nu=1/2$ \cite{GaussProc}, 
the Mat\'{e}rn covariance function reduces to the exponential covariance model and describes a rough field.
The value $\nu=\infty$ corresponds to a Gaussian covariance model, which describes a very 
smooth field, that is infinitely differentiable.
Random fields with a Mat\'{e}rn covariance function are $\lfloor \nu-1 \rfloor$ times mean square differentiable.

Thus, the
hyperparameter $\nu$ controls the degree of smoothness. 

The $d$-dimensional Fourier transform $\F^{d}(C(r,\nu)$ of the Mat\'{e}rn  
kernel, defined in \refeq{eq:Matern}, in $\mathbb{R}^d$ is given by \cite{Matern1986a}
\begin{equation}
\label{eq:FFTMatern}
\U(\xi):=\F^{d}(C(r,\nu)=\beta \cdot \left( 1+\frac{\ell^2}{2\nu}\vert \xi \vert^2\right)^{-\nu-d/2},
\end{equation}
where $\beta=\beta(\nu, \ell,n)$ is a constant 
and $\vert \xi \vert$ is the Euclidean distance in $\mathbb{R}^d$. 

The following tensor approach also applies to the case of anisotropic distance, where
$r^2= \sqrt{\langle A(x-y),(x-y)\rangle}$, and $A$ is a positive diagonal $d\times d$ matrix.

\section{Low-rank tensor decompositions}
\label{sec:low-rank}
In this section, we review the definitions of the CP and Tucker tensor formats.
Then, we provide the analytic $\sinc$-based 
proof of the existence of low-rank tensor approximations of Mat\'{e}rn functions.
We investigate numerically the behavior of
the Tucker and CP ranks across a wide range of parameters specific to the family of Mat\'{e}rn 
kernels in Eq.~\ref{eq:FFTMatern}. The Tucker tensor format is used in this work for additional 
rank compression of the CP factors. 
There are no reliable algorithms to compute CP decomposition, which can be difficult to compute, 
but there are such algorithms 
for Tucker decomposition. The Tucker decomposition is only limited  w.r.t. the available memory storage, 
since the term $r^d$ 
in $\mathcal{O}(drn +r^d)$ grows exponentially with $d$.
\subsection{General definitions}
 \label{ssec:low-rank_generaldef}
CP and Tucker rank-structured tensor formats have been 
applied for the quantitative analysis of correlation in multi-dimensional experimental data
for a long time in
chemometrics and signal processing \cite{smilde-book-2004,Cichocki:2002}. 
The Tucker tensor format was introduced in 1966 for tensor decomposition of
multidimensional arrays in chemometrics \cite{Tuck:66}. 
Though the canonical representation of multivariate functions was
introduced as  early as in 1927 \cite{Hitch:27}, only the Tucker tensor format provides a
stable algorithm for decomposition of full-size tensors. A mathematical approval 
of the Tucker decomposition algorithm was presented in papers on
higher order singular value decomposition (HOSVD) and the Tucker
  ALS algorithm for orthogonal Tucker approximation of higher 
order tensors  \cite{DMV-SIAM2:00}. 
For higher dimensions, the so-called Matrix Product States (MPS) 
(see the survey paper \cite{Scholl:11})
or the Tensor Train (TT) \cite{Osel_TT:11} decompositions can be applied. 
However, for 3D applications, the Tucker and CP
tensor formats remain the best choices.
The fast convergence of the Tucker decomposition was proved and demonstrated numerically for 
higher order tensors that 
arise from the discretization of linear operators and functions in $\mathbb{R}^d$ for a class 
of function-related tensors 
and Green's kernels in particular, it was found that the approximation error of the Tucker 
decomposition decayed 
exponentially in the Tucker rank \cite{khor-rstruct-2006,Khoromskij_Low_Tacker}.

These results inspired the canonical-to-Tucker (C2T)
and Tucker-to-canonical (T2C) decompositions 
for function-related tensors in the case of large input ranks, as well as the multigrid 
Tucker approximation \cite{khor-ml-2009}.

A tensor of order $d$ in a full format is defined as a multidimensional array
over a $d$-tuple index set:
\begin{equation}
 \label{Tensor_def}
{\bf A}=[a_{i_1,\ldots,i_d}]\equiv [a(i_1,\ldots,i_d)] \; \in 
\mathbb{R}^{n_1 \times \ldots \times n_d}\;
\mbox{  with } \quad i_\ell\in I_\ell:=\{1,\ldots,n_\ell\}.
\end{equation}
Here,
$\bf A$ is an element of the linear space
\begin{equation*}
 \mathbb{V}_n=\bigotimes_{\ell=1}^d \mathbb{V}_{\ell},\quad \mathbb{V}_{\ell}=\mathbb{R}^{I_\ell}
\end{equation*}
equipped with the Euclidean scalar product $\langle \cdot,\cdot \rangle:
\mathbb{V}_{n}\times \mathbb{V}_{n}\rightarrow \mathbb{R}$, defined as
\begin{equation*}
 \langle {\bf A},{\bf B} \rangle:=\sum_{(i_1...i_d)\in \mathcal{I}}a_{i_1...i_d}b_{i_1...i_d},
 \quad \text{for}\, \bf{A},\,\bf{B}\in \mathbb{V}_n.
\end{equation*}

Tensors with all dimensions having equal size $n_\ell =n$, $\ell=1,\ldots d$,
are called $n^{\otimes d}$ tensors. The  required storage size 
scales exponentially with the dimension, $n^d$, which results in the 
so-called ``curse of dimensionality''. 

To avoid exponential scaling in the dimension, the 
rank-structured separable representations (approximations) of the multidimensional tensors can be used.
The simplest separable element is given by the rank-$1$ tensor, 
\[ {\bf U} = {\bf u}^{(1)}\otimes \ldots 
 \otimes {\bf u}^{(d)}\in \mathbb{R}^{n_1 \times \ldots \times n_d},
\]
with entries
$
 u_{i_1,\ldots, i_d}=  u^{(1)}_{i_1} \cdots u^{(d)}_{i_d},
$
which requires only $n_1+ \ldots +n_d$ numbers for storage. 

The rank-1 canonical tensor is a discrete counterpart of the separable $d$-variate function,
which can be represented as the product of univariate functions,
\[
 f(x_1,x_2,\ldots , x_d)= f_1(x_1) f_2(x_2)\ldots f_d(x_d).
\]
An example of the separable $d$-variate function is $f(x_1,x_2,x_3) =e^{(x_1 + x_2+ x_3) }$.
Then, by discretization of this multivariate function on a tensor grid in a computational box, we obtain
a canonical rank-1 tensor.

\begin{figure}[htb]
\centering
\includegraphics[width=6.5cm]{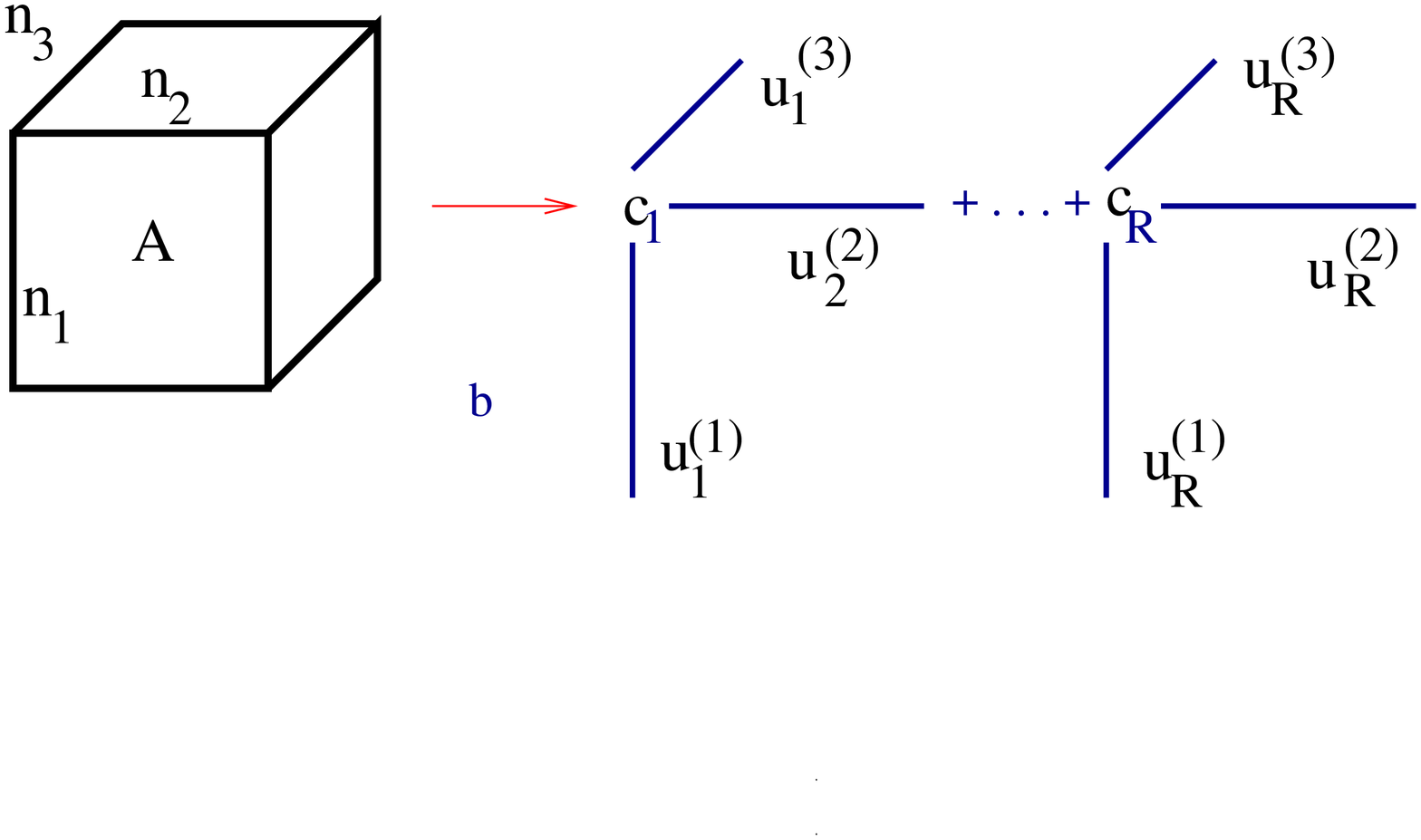}\quad\quad\quad
\includegraphics[width=6.5cm]{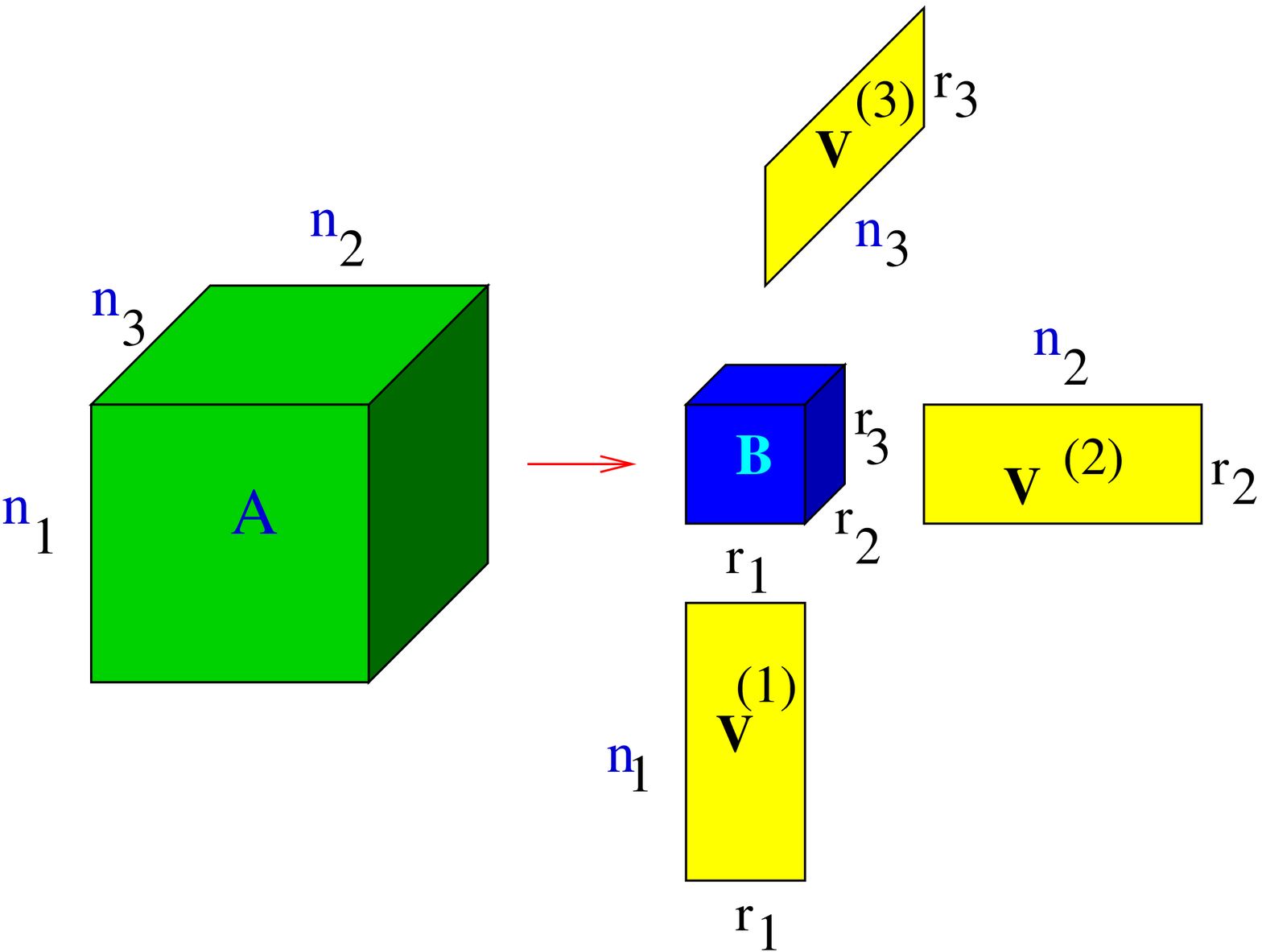}
\caption{\small Canonical (left) and Tucker (right) decompositions of 3D tensors.}
\label{fig:CP_Tucker}
\end{figure}

A tensor in the $R$-term canonical format is defined by a finite sum of rank-$1$ tensors 
(Fig.~\ref{fig:CP_Tucker}, left),
 \begin{equation}\label{eqn:CP_form}
   {\bf A}_{c} = {\sum}_{k =1}^{R} \xi_k
   {\bf u}_k^{(1)}  \otimes \ldots \otimes {\bf u}_k^{(d)},  \quad  \xi_k \in \mathbb{R},
\end{equation}
where ${\bf u}_k^{(\ell)}\in \mathbb{R}^{n_\ell}$ are normalized vectors, 
and $R$ is the canonical rank. The storage cost of this
parametrization is bounded by  $d R n$.
An element  $a(i_1,...,i_d)$ of tensor ${\bf A}=\sum_{i=1}^R\bigotimes_{\nu=1}^d u_{i \nu}$ 
can be computed as
\begin{equation*}
 a(i_1,...,i_d)= \sum_{\alpha=1}^R u_1(i_1,\alpha)u_2(i_2,\alpha)...u_d(i_d,\alpha).
\end{equation*}

An alternative (contracted product) notation is used in computer science community:
\begin{equation}
 \label{Can_Contr}
 {\bf A} = {\bf C} \times_1 U^{(1)} \times_2 U^{(2)} \times_3 \cdots \times_d U^{(d)},
 \end{equation}
where ${\bf C}=diag\{c_1,...,c_d\} \in\mathbb{R}^{R^{\otimes d}}$, and
$ U^{(\ell)}=[{\bf u}_1^{(\ell)}\ldots {\bf u}_R^{(\ell)} ]\in\mathbb{R}^{n_\ell \times R}$.
An analogous multivariate function can be represented by a sum
of univariate functions:
\[
 f(x_1,x_2,\ldots , x_d) = \sum\limits^R_{k=1} f_{1,k}(x_1) f_{2,k}(x_2)\ldots f_{d,k}(x_2).
\]
For $d\geq 3$, there are no stable algorithms to compute the 
canonical rank of a tensor ${\bf A}$, that is
the minimal number $R$ in representation (\ref{eqn:CP_form}), and the respective
decomposition with the polynomial cost in $d$, i.e., the computation of the
canonical decomposition is an $N$-$P$ hard problem \cite{Hastad:90}.

The Tucker tensor format (Fig.~\ref{fig:CP_Tucker}, right) is suitable for stable numerical 
decompositions with a fixed
truncation threshold.
We say that the tensor ${\bf A} $ is represented in the rank-$\bf r$ orthogonal Tucker format 
with the rank parameter ${\bf r}=(r_1,\ldots,r_d)$ if 
\begin{equation}\label{eqn:Tucker_form}
  {\bf A}  =\sum\limits_{\nu_1 =1}^{r_1}\ldots
\sum\limits^{r_d}_{{\nu_d}=1} \beta_{\nu_1, \ldots ,\nu_d}
\,  {\bf v}^{(1)}_{\nu_1} \otimes \ldots \otimes {\bf v}^{(\ell)}_{\nu_{\ell}} 
\ldots \otimes {\bf v}^{(d)}_{\nu_d},\quad \ell=1,\ldots,d, 
\end{equation}
where $\{{\bf v}^{(\ell)}_{\nu_\ell}\}_{\nu_\ell=1}^{r_\ell}\in \mathbb{R}^{n_\ell}$ 
represents a set of orthonormal vectors for $\ell=1,\ldots,d$, 
and $\bB=[\bB_{\nu_1,\ldots,\nu_d}] \in \mathbb{R}^{r_1\times \cdots \times r_d}$ is 
the Tucker core tensor. 
The storage cost for the Tucker tensor is bounded by $d r n +r^d$, 
with $r=|{\bf r}|:=\max_\ell r_\ell$.
Using the orthogonal side matrices $V^{(\ell)}=[v^{(\ell)}_1 \ldots v^{(\ell)}_{r_\ell}]$ 
and contracted products, the Tucker tensor decomposition may be presented in the alternative notation,
\[
  { {\bf A}_{({\bf r})}=\bB \times_1 V^{(1)}\times_2 V^{(2)} \times_3 \ldots 
   \times_d V^{(d)}.} 
\]  
In the case $d=2$, the orthogonal Tucker decomposition is equivalent to the singular 
value decomposition (SVD) of a rectangular matrix. 

\subsection{Tucker decomposition of full format tensors}
\label{ssec:Tucker_decomp}
 
\begin{figure}[htb]
\centering
\includegraphics[width=10.0cm]{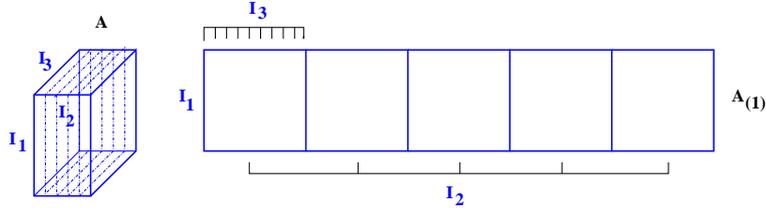} 
\caption{\small Unfolding of a 3D tensor along the mode $I_\ell$ with $\ell=1$.}
\label{fig:Unfolding}
\end{figure}

We use the following algorithm to compute the Tucker decomposition of the full format tensor. 
The most time-consuming 
part of the Tucker algorithm is higher order singular value decomposition (HOSVD), the computation of 
the initial guess for matrices $V^{(\ell)}$ using
the SVD of the matrix unfolding $A_{(\ell)}$, $\ell=1,\; 2, \; 3 $ 
(Fig.~\ref{fig:Unfolding}), of the original tensor 
along each mode of a  tensor \cite{DMV-SIAM2:00}. 
Figure \ref{fig:Unfolding} illustrates the matrix unfolding of the full format 
tensor ${\bf A}$ (\ref{Tensor_def}) along the index set $I_1=\{1,\ldots,n_1\}$. 
 
The second part of the algorithm
is the ALS procedure. For every tensor mode, a ``single-hole'' tensor of reduced size
is constructed by the mapping all of the modes 
of the original tensor except one into the subspaces $V^{(\ell)}$. 
Then, the subspace $V^{(\ell)}$ for the current mode is updated by SVD
of the unfolding of the ``single hole'' tensor for this mode.
This alternates over all modes of the tensor, which are updated
at the current iteration of ALS. The final step of the algorithm is computation of the core tensor by using the
ultimate mapping matrices from ALS.
    
 The numerical cost of Tucker decomposition for full size tensors is dominated by
 the initial guess, which is estimated as
 $O(n^{d+1})$ when all $n_\ell=n$ are equal , or $O(n^{4})$ for our 3D case. 
 This step restricts the available size of the tensor to be decomposed since, 
 for conventional computers, the 3D case $n_\ell> 10^2$ is the limiting case for SVD.
 
The multigrid Tucker algorithm for full size tensors allows the computational complexity to be 
linear in the full size of the tensor, $O(n^d)$, i.e. $O(n^3)$ for 3D tensors \cite{khor-ml-2009}. 
It is computed on a sequence of diadically refined grids and
is based on implementing the HOSVD only at the coarsest grid level, $n_0 \ll n$.
 The initial guess for the ALS procedure is computed at each refined 
level by the interpolation of the dominating Tucker subspaces 
obtained from the previous coarser grid. In this way, at fine 3D Cartesian grids, we need 
29only  $O(n^3)$ storage (to represent the initial full format tensor) 
to contract with the Tucker side matrices, obtained by 
the Tucker approximation via ALS on the previous grids.
\subsection{Illustration of the low-rank approximation idea}
\label{sec:idea}
Here, we describe a possible ways to find a low-rank tensor approximation of the 
Mat\'{e}rn covariance matrix $C(r,\nu)$ (Diagram 1).
We let $\F^{d}=\bigotimes_{\nu=1}^d \F_{\nu}$ be the $d$-dimensional Fourier transform, where
$\F^{-d}= \bigotimes_{\nu=1}^d \F^{-1}_{\nu} $ is its inverse and $\bigotimes$ denotes the Kronecker product.
We assume that $\U(\xi)=\F^{d}(C(r,\nu))$ is known analytically and has a low-rank tensor 
approximation $\mathbf{U}= \sum_{j=1}^{r} \bigotimes_{\nu=1}^d \mathbf{u}_{j \nu}$.

Since the
Fourier and inverse Fourier Transformations do not change the Kronecker tensor 
rank of the argument \cite{nowak2013kriging},
then by applying the inverse Fourier, we obtain a low-rank representation of the 
covariance function by applying the inverse Fourier:
\begin{equation}
\F^{-d}\left(\mathbf{U}\right)
=\left( \bigotimes_{\nu=1}^d \F^{-1}_{\nu} \right) \sum_{i=1}^{r} 
\left( \bigotimes_{\nu=1}^d \mathbf{u}_{\nu i} \right)
=\sum_{i=1}^{r}\bigotimes_{\nu=1}^d\left(\F^{-1}_{\nu} \left(\mathbf{u}_{\nu i}\right)\right)=
\sum_{i=1}^{r}\bigotimes_{\nu=1}^d\tilde{\mathbf{u}}_{\nu i}=:C(r,\nu),
\end{equation}
where $\tilde{\mathbf{u}}_{\nu i}:=\F^{-1}_{\nu} \left(\mathbf{u}_{\nu i}\right)$.

\subsection{$\mbox{Sinc}$ approximation of the Mat\'ern function}
\label{ssect: Sinc_C2T}
 {The $\mbox{Sinc}$ method provides a  constructive approximation 
of the multivariate functions in the form of a low-rank canonical representation.
It can be also used for the theoretical proof and for the rank estimation.}
Methods for the separable approximation of the 3D Newton kernel and many other spherically 
symmetric functions that use the Gaussian sums have been developed 
since the initial studies in chemical \cite{Boys:56} and mathematical 
literature \cite{Stenger,Braess:BookApTh:86,HaKhtens:04I,GHK:05}.
Here, we use a tensor-decomposition approach for lattice-structured interaction 
potentials \cite{VeKhor_NLLA:15}.
We also recall the grid-based method for a low-rank canonical  
representation of the spherically symmetric kernel function $q(\|x\|)$, 
where $x\in \mathbb{R}^d$ where $d=2,3,\ldots$, by its projection onto the set
of piecewise constant basis functions; see \cite{BeHaKh:08} for the case of 
Newton and Yukawa kernels for $x\in \mathbb{R}^3$.

Following the standard schemes, we introduce 
the uniform $n \times n \times n$ rectangular 
Cartesian grid $\Omega_{n}$ with mesh size $h=2b/n$ (we assume that $n$ is even) 
in the computational domain $\Omega=[-b,b]^3$.
We let $\{ \psi_\textbf{i}\}$ be a set of tensor-product piecewise constant basis functions,
$ \psi_\textbf{i}(\textbf{x})=\prod_{\ell=1}^3 \psi_{i_\ell}^{(\ell)}(x_\ell)$,
for the $3$-tuple index ${\bf i}=(i_1,i_2,i_3)\in {\cal I}$, ${\cal I}=I_1\times I_2 \times I_3$, with
$i_\ell \in I_\ell=\{1,...,n\}$, where $\ell=1,\, 2,\, 3 $.
The generating kernel $q(\|x\|)$ is discretized by its projection onto the basis 
set $\{ \psi_\textbf{i}\}$
in the form of a third-order tensor of size $n\times n \times n$, which is defined entry-wise as
\begin{equation}  \label{galten}
\mathbf{Q}:=[q_{\bf i}] \in \mathbb{R}^{n\times n \times n},  \quad
 q_{\bf i} = 
\int_{\mathbb{R}^3} \psi_{\bf i} ({x}) q(\|{x}\|) \,\, \mathrm{d}{x}.
\end{equation}

The low-rank canonical decomposition of the third-order tensor $\mathbf{Q}$ is based 
on applying exponentially convergent 
$\operatorname*{\sinc}$-quadratures to the integral representation of  
the function $q(p)$, $p \in \mathbb{R}$ in the form
\[
 q(p)= \int_{\mathbb{R}} a_1(t) e^{- p^2 a_2(t)} dt,
\]
specified by the weights $a_1(t), a_2(t)  >0$.
Diagram 2 illustrates a scheme of the proof of existence of the canonical low-rank tensor approximation. 
It could be easier to apply the Laplace transform to the Fourier transform of a Mat\'{e}rn 
covariance matrix, than to the Mat\'{e}rn covariance. 
To approximate the resulting Laplace integral we apply the $\sinc$ quadrature. The number $(2M+1)$ of terms 
in the approximate sum (Eq.~\ref{eqn:laplace}) is the canonical tensor rank.

\begin{tikzpicture}
\node[shortblock](process00){$\approx \sum_{i=1}^r \bigotimes_{\mu=1}^d f_{i \mu} $};
\node[evshortblock, left of=process00,xshift=-5em, yshift=-0em](process0){$\mathcal{L}(s)$};
\node[evshortblock, left of=process0, xshift=-2em, yshift=-0em]
(process2){$f_{\alpha,\nu}(\rho)$};
\node[evshortblock, left of=process2, xshift=-3em, yshift=-0em]
(process3){$C_{\nu,\ell}(r)$};
\path[line2](process0)--node[yshift=1em]{$\sinc$}(process00);
\path[line2](process2)--node[yshift=1em]{Laplace }(process0);
\path[line2](process3)--node[yshift=1em]{Fourier}(process2);
\node [below=1cm, xshift=17em, align=flush center,text width=15cm] at (process3)
{Diagram 2: Scheme of the proof of existence of low-rank tensor approximation, $r=2M+1$.};
\end{tikzpicture}

In particular,  the $\operatorname*{\sinc}$-quadrature for the Laplace-Gauss transform, 
\begin{align} \label{eqn:laplace} 
q(p)=\int_{\mathbb{R}_+} a(t) e^{- t^2 p^2} \,\mathrm{d}t \approx
\sum_{k=-M}^{M} a_k e^{- t_k^2 p^2} \quad \mbox{for} \quad |p| > 0,\quad p \in \mathbb{R},
\end{align} 
can be applied, where the quadrature points $(t_k)$ and weights $(a_k)$ are given by 
\begin{equation} \label{eqn:hM}
t_k=k \mathfrak{h}_M, \;\;\text{and} \quad a_k=a(t_k) \mathfrak{h}_M, \quad \text{where}\;\;
\mathfrak{h}_M=C_0 \log(M)/M , \quad C_0>0.
\end{equation}
Under the assumption $0< a_0 \leq |p |  < \infty$,
this quadrature can be proven to provide an exponential convergence rate in $M$ (uniformly in $p$)
for a class of functions $a(z)$, that are analytic
in a certain strip $|z|\leq D $ of the complex plane, such that the functions
$a_1(t) e^{- p^2 a_2(t)}$ decay polynomially or exponentially on the real axis.  
The exponential convergence of the $\sinc$-approximation in the number of terms 
(i.e., the canonical rank $R=2M+1$) was analyzed elsewhere \cite{Stenger,Braess:BookApTh:86,HaKhtens:04I}.

We assume that a representation similar to (\ref{eqn:laplace}) exists 
for any fixed $x=(x_1,x_2,x_3)\in \mathbb{R}^3$, 
such that $\|{x}\| > a_0 > 0$. 
Then, we apply the $\operatorname*{\sinc}$-quadrature approximation (\ref{eqn:laplace}) and (\ref{eqn:hM}) 
to obtain the separable expansion
\begin{equation} \label{eqn:sinc_Newt}
 q({\|{x}\|}) =   \int_{\mathbb{R}_+} a(t)
e^{- t^2\|{x}\|^2} \,\mathrm{d}t  \approx 
\sum_{k=-M}^{M} a_k e^{- t_k^2\|{x}\|^2}= 
\sum_{k=-M}^{M} a_k  \prod_{\ell=1}^3 e^{-t_k^2 x_\ell^2},
\end{equation}
providing an exponential convergence rate in $M$:
\begin{equation} \label{sinc_conv}
\left|q({\|{x}\|}) - \sum_{k=-M}^{M} a_k e^{- t_k^2\|{x}\|^2} \right|  
\le \frac{C}{a}\, \displaystyle{e}^{-\beta \sqrt{M}},  
\quad \text{with some} \ C,\beta >0.
\end{equation}
By combining \eqref{galten} and \eqref{eqn:sinc_Newt}, and taking into account the 
separability of the Gaussian basis functions, we arrive at the low-rank 
approximation of each entry of the tensor $\mathbf{Q}$:
\begin{equation*} \label{eqn:C_nD_0}
 q_{\bf i} \approx \sum_{k=-M}^{M} a_k   \int_{\mathbb{R}^3}
 \psi_{\bf i}({x}) e^{- t_k^2\|{x}\|^2} \mathrm{d}{x}
=  \sum_{k=-M}^{M} a_k  \prod_{\ell=1}^{3}  \int_{\mathbb{R}}
\psi^{(\ell)}_{i_\ell}(x_\ell) e^{- t_k^2 x^2_\ell } \mathrm{d} x_\ell.
\end{equation*}

Recalling that $a_k > 0$, we define the vector $\textbf{q}$ as
\begin{equation*} \label{eqn:galten_int}
\textbf{q}^{(\ell)}_k
= a_k^{1/3} \left[b^{(\ell)}_{i_\ell}(t_k)\right]_{i_\ell=1}^{n_\ell} \in \mathbb{R}^{n_\ell}
\quad \text{with } \quad b^{(\ell)}_{i_\ell}(t_k)= 
\int_{\mathbb{R}} \psi^{(\ell)}_{i_\ell}(x_\ell) e^{- t_k^2 x^2_\ell } \mathrm{d}x_\ell.
\end{equation*}
Then, the third order tensor $\mathbf{Q}$ can be approximated by 
the $R$-term ($R=2M+1$) canonical representation
\begin{equation} \label{eqn:sinc_general}
    \mathbf{Q} \approx  \mathbf{Q}_R =
\sum_{k=-M}^{M} a_k \bigotimes_{\ell=1}^{3}  {\bf b}^{(\ell)}(t_k)
= \sum\limits_{k=-M}^{M} {\bf q}^{(1)}_k \otimes {\bf q}^{(2)}_k \otimes {\bf q}^{(3)}_k
\in \mathbb{R}^{n\times n \times n}, 
\end{equation}
where ${\bf q}^{(\ell)}_k \in \mathbb{R}^n$. 
Given a threshold $\varepsilon >0 $,  $M$ can be chosen as the minimal number
such that in the max-norm
\begin{equation*} \label{eqn:error_control}
\| \mathbf{Q} - \mathbf{Q}_R \|  \le \varepsilon \| \mathbf{Q}\|.
\end{equation*}
The skeleton vectors can be reindex by $k \mapsto k'=k+M+1$, 
${\bf q}^{(\ell)}_k \mapsto {\bf q}^{(\ell)}_{k'}$, ($k'=1,...,R=2M+1$), $\ell=1, 2, 3$.
The symmetric canonical tensor, ${\bf Q}_{R}\in \mathbb{R}^{n\times n\times n}$ in (\ref{eqn:sinc_general}),
approximates the 3D symmetric kernel function 
$q({\|x\|})$ ($x\in \Omega$), centered at the origin, such that 
${\bf q}^{(1)}_{k'}={\bf q}^{(2)}_{k'}={\bf q}^{(3)}_{k'}$ ($k'=1,...,R$).

In some applications, the tensor can be given in the canonical tensor format, but with large
rank $R$ and discretized on large grids $n\times n \times \ldots \times n$; thus,
computation of the initial of guess in the Tucker-ALS decomposition algorithm becomes intractable.
This situation may arise when composing the tensor approximation of complicated kernel
functions from simple radial functions that can be represented 
in the low-rank CP format.

For such cases, the canonical-to-Tucker decomposition algorithm was introduced \cite{khor-ml-2009}.
It is based on the minimization ALS procedure, similar to the Tucker algorithm, 
described in Section \ref{ssec:Tucker_decomp} for full-size tensors, but the 
initial guess is computed by just the 
SVD of the side matrices, 
$ U^{(\ell)}=[{\bf u}_1^{(\ell)}\ldots {\bf u}_R^{(\ell)} ]\in\mathbb{R}^{n_\ell \times R}$,
 $\ell=1,\; 2,\; 3$; see (\ref{Can_Contr}). This schema is the Reduced HOSVD
(RHOSVD), which does not require unfolding of the
full tensor.
 
Another efficient rank-structured representation of the multidimensional tensors 
is the mixed-tensor format \cite{VeKh_Diss:10}, which combines either the
canonical-to-Tucker decomposition with the Tucker-to-canonical decomposition, 
or standard Tucker decomposition with the canonical-to-Tucker decomposition, in order
to produce a canonical tensor from a full-size tensor.

\subsection{Laplace transform of the covariance matrix}
\label{ssec:Laplace_transf}
The integral representations like (\ref{eqn:laplace}) can be derived by the Laplace 
transform  either directly to the Mat{\'e}rn covariance function or to its spectral density 
 (\ref{eq:FFTMatern}).
 
 For example, in the case of the Newton kernel, $q(p)=1/p$, and,the Laplace-Gauss transform representation 
takes the form
\begin{equation}
 \frac{1}{p}= \frac{2}{\sqrt{\pi}}\int_{\mathbb{R}_+} e^{- p^2 t^2 } dt, \quad 
 \mbox{where}\quad p=\|x\|=\sqrt{x_1^2 + x_2^2  + x_3^2}.
\end{equation}
In this case, ${\bf q}^{(\ell)}_k={\bf q}^{(\ell)}_{-k}$, 
and the sum of (\ref{eqn:sinc_general}) reduces to $k=0,1,...,M$, implying that $R=M+1$.
Therefore, the Laplace transform representation of the Slater function $q(p)=e^{-2\sqrt{\alpha p}}$ 
(i.e., exponential covariance)
with $p=\|x\|^2$ can be written as
\begin{equation}
  q(p)=e^{-2\sqrt{\alpha p}} =\frac{\sqrt{\alpha}}{\sqrt{\pi}}
  \int_{\mathbb{R}_{+}} t^{-3/2} exp(-\alpha/t - p t) d t.
\end{equation}

When the Mat{\'e}rn spectral density in (\ref{eq:FFTMatern}) has an
even dimension parameter $d=2d_1$, $d_1=1,2,...$ and $\nu=0,1,2,...$, the Laplace transform,
\begin{equation}
 \frac{\eta!}{(p+a)^{\eta+1}}=\int_{\mathbb{R}_{+}} t^{\eta} e^{-at} e^{-p t} dt,
\end{equation}
can be applied after substituting $p= \|\xi\|^2$ in $q(p)=\beta \left 
( 1+\frac{\ell^2}{2\nu}p\right)^{\eta}$, $\eta=-\nu-d_1$. 

If $-\eta=\nu+d_1=1/2,3/2,5/2,...,\frac{2k+1}{2},...$, then the Laplace transform is
\begin{equation}
 \frac{(2\eta)! \sqrt{\pi}}{\eta! 4^{\eta}} \frac{1}{(p+a)^{\eta}\sqrt{p+a}}=
 \int_{\mathbb{R}_{+}} \frac{t^{\eta}}{\sqrt{t}} e^{-at} e^{-p t} dt, \quad \eta\in \mathbb{N}.
\end{equation}

%
%
\subsection{Covariance matrix in rank-structured tensor format}
\label{ssec:Cover_tensor}

  In what follows, we consider the CP approximation of the radial function $q(r)$ in the positive sector,
  i.e., on the domain $[0,b]^3$ (by symmetry, the canonical tensor can be extended 
  to the whole computational domain $[-b, b]^3$.)
 We let the covariance function $q(r)=C_{\nu,\ell}(r)$ in (\ref{eq:Matern}) be 
 represented by the rank-$R$ symmetric CP tensor on an $n\times n \times n$ tensor grid, 
 denoted by $\Omega_n\subset \Omega=[0,b]^3$  as described  in the previous sections
 \begin{equation} \label{eqn:Covar_CP}
  q(r) \mapsto  \mathbb{Q} \approx  \mathbb{Q}_R =
\sum\limits_{k=1}^{R} {\bf q}^{(1)}_k \otimes {\bf q}^{(2)}_k \otimes {\bf q}^{(3)}_k
\in \mathbb{R}^{n\times n \times n}, 
\end{equation}
with the same skeleton vectors ${\bf q}^{(\ell)}_k \in \mathbb{R}^n$ for $\ell=1,2,3$. 
 
We define the covariance matrix ${\bf C}=[c_{{\bf i},{\bf j}}]\in \mathbb{R}^{{\bf n}\times {\bf n}}$ entry-wise by 
 \begin{equation}
  c_{{\bf i},{\bf j}}= C_{\nu,\ell}(\|x_{\bf i} - y_{\bf j} \|), \quad {\bf i},{\bf j}\in {\cal I}.
\end{equation}
Using the tensor representation (\ref{eqn:Covar_CP}), we represent the large $n^3\times n^3$
matrix ${\bf C}$ in the rank-$R$ Kronecker (tensor) format as
 \begin{equation} \label{eqn:CovarMatr_CP}
 {\bf C} \approx {\bf {C}}_R= 
 \sum\limits_{k=1}^{R} {\bf Q}^{(1)}_k \otimes {\bf Q}^{(2)}_k \otimes {\bf Q}^{(3)}_k,
\end{equation}
where the symmetric Toeplitz matrix, 
${\bf Q}^{(\ell)}_k=\mbox{Toepl}[{\bf q}^{(\ell)}_k]\in \mathbb{R}^{n\times n}$,
$\ell=1,2,3$, is defined by its first column, which is specified by the skeleton 
vectors ${\bf q}^{(1)}_k={\bf q}^{(2)}_k={\bf q}^{(3)}_k$ in the decomposition (\ref{eqn:Covar_CP}).

\begin{figure}[htb]
\centering
\includegraphics[width=7.5cm]{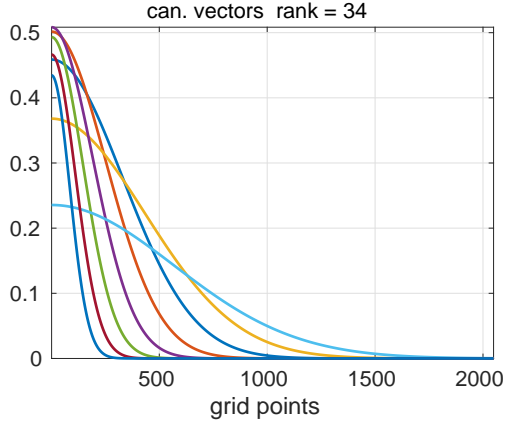}
\caption{\small Selected eight canonical vectors from the full 
set ${\bf q}^{(1)}_k$, $k=1,\ldots,R$, see (\ref{eqn:Covar_CP}).}
\label{fig:gener_funct}
\end{figure}

Figure \ref{fig:gener_funct} illustrates eight selected canonical generating vectors from 
${\bf q}^{(1)}_k$ for $k=1,\ldots,R$, $R=34$, on a grid of size $n=2049$  
for the Slater function $e^{-\|x\|^p}$, which defines the corresponding Toeplitz matrices.

A Toeplitz matrix can be multiplied by a vector in $O(n \log n)$ operations via 
complementing it to the circulant matrix. In general, the inverse of a Toeplitz matrix 
cannot be calculated in the closed form; unlike circulant matrices,
which can be diagonalized by the Fourier transform.

In the rest of this section, we introduce the numerical scheme, based on certain
specific properties of the skeleton matrices in the symmetric rank-$R$ decomposition 
(\ref{eqn:CovarMatr_CP}) that is capable of 
rank-structured calculations
of the analytic matrix functions ${\cal F}({\bf C}_R)$. 
We discuss the most interesting examples of the functions 
${\cal F}_1({\bf C})={\bf C}^{-1}$ and ${\cal F}_2({\bf C})={\bf C}^{1/2}$, 
where ${\bf C}={\bf C}_R$.

Given an symmetric positive definite matrix such that $\|{\bf A}\|=q<1$,  
the matrix-valued function is given as the exponentially fast converging series   
\begin{equation}
\label{eq:cF}
 {\cal F}({\bf A})= {\bf E} + a_1 {\bf A} +a_2 {\bf A}^2 + \ldots,
\end{equation}
where the matrix ${\bf A}$, acting in the multi-dimensional index set, allows the low-rank
Kronecker tensor decomposition. Then, the low-rank tensor approximation of
${\cal F}({\bf A})$ can be computed by the ``add-and-compress" scheme, 
where each term in the series above (\ref{eq:cF}) is summed using the rank-truncation 
algorithm in the corresponding format.

To limit the rank-structured evaluation of ${\cal F}_1({\bf C})$ to the described framework,
we propose the special rank-structured additive splitting of the covariance matrix ${\bf C}$
with the easily invertible dominating part. To that end, we construct the diagonal matrix
${\bf Q}^{(1)}_0$ by assembling all of the diagonal sub-matrices in ${\bf Q}^{(1)}_k$, $k=1,\ldots,R$, 
(in the following, we simplify the notation by omitting the upper index $(1)$):
\begin{equation}
 {\bf Q}_0 := \sum\limits_{k=1}^R \mbox{diag} ({\bf Q}_k),
\end{equation}
and modify each  Toeplitz matrix ${\bf Q}_k$ by subtracting its diagonal part,
\begin{equation}
 {\bf Q}_k \mapsto \widehat{\bf Q}_k:= {\bf Q}_k - \mbox{diag} ({\bf Q}_k), \quad k=1,\ldots,R.
\end{equation}
Using the matrices defined above, we introduce the rank $R+1$ additive splitting of ${\bf C}$ , 
which is defined by
the skeleton matrices ${\bf Q}_0$ and $\widehat{\bf Q}_k$ since ${\bf {C}}$ is Kronecker symmetrical,
\begin{equation} \label{eqn:CovarMatr_CP_split}
 {\bf {C}}= {\bf Q}_0 \otimes {\bf Q}_0 \otimes {\bf Q}_0 + 
 \sum\limits_{k=1}^{R} \widehat{\bf Q}^{(1)}_k \otimes \widehat{\bf Q}^{(2)}_k 
 \otimes \widehat{\bf Q}^{(3)}_k.
\end{equation}
Hence, we have
\begin{equation}
 {\bf {C}}^{-1}={\bf Q}_0^{-1} \otimes {\bf Q}_0^{-1} \otimes {\bf Q}_0^{-1}({\bf {E}} + \sum\limits_{k=1}^{R} 
{\bf Q}_0^{-1} \widehat{\bf Q}^{(1)}_k \otimes {\bf Q}_0^{-1}\widehat{\bf Q}^{(2)}_k 
\otimes {\bf Q}_0^{-1}\widehat{\bf Q}^{(3)}_k)^{-1}.
\end{equation}
Likewise, since ${\bf Q}_0 \otimes {\bf Q}_0 \otimes {\bf Q}_0$ is a scaled identity, we obtain
\begin{equation}
\label{eq:sqrt}
 {\bf {C}}^{1/2}={\bf Q}_0^{1/2} \otimes {\bf Q}_0^{1/2} \otimes {\bf Q}_0^{1/2}
 ({\bf {E}} + \sum\limits_{k=1}^{R} 
{\bf Q}_0^{-1} \widehat{\bf Q}^{(1)}_k \otimes {\bf Q}_0^{-1}\widehat{\bf Q}^{(2)}_k 
\otimes {\bf Q}_0^{-1}\widehat{\bf Q}^{(3)}_k)^{1/2}.
\end{equation}
We assume that $\|{\bf Q}_0^{-1} \widehat{\bf Q}^{(1)}_k \|<1$ for  $k=1,\ldots,R$ in 
some norm; thus, we can apply the ``add-and-compress" scheme described above.

\begin{figure}[htb]
\centering
\includegraphics[width=7.5cm]{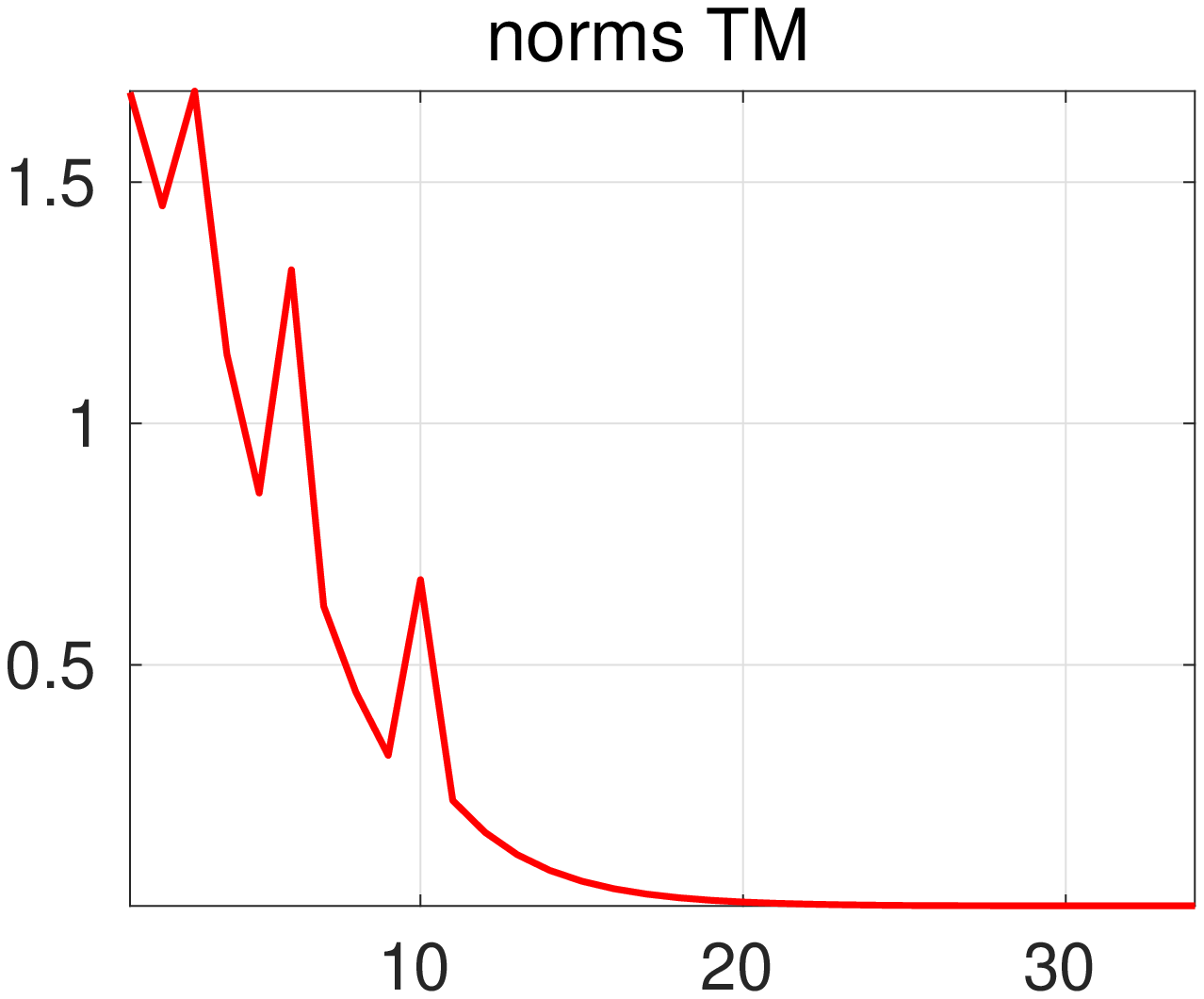}
\includegraphics[width=7.5cm]{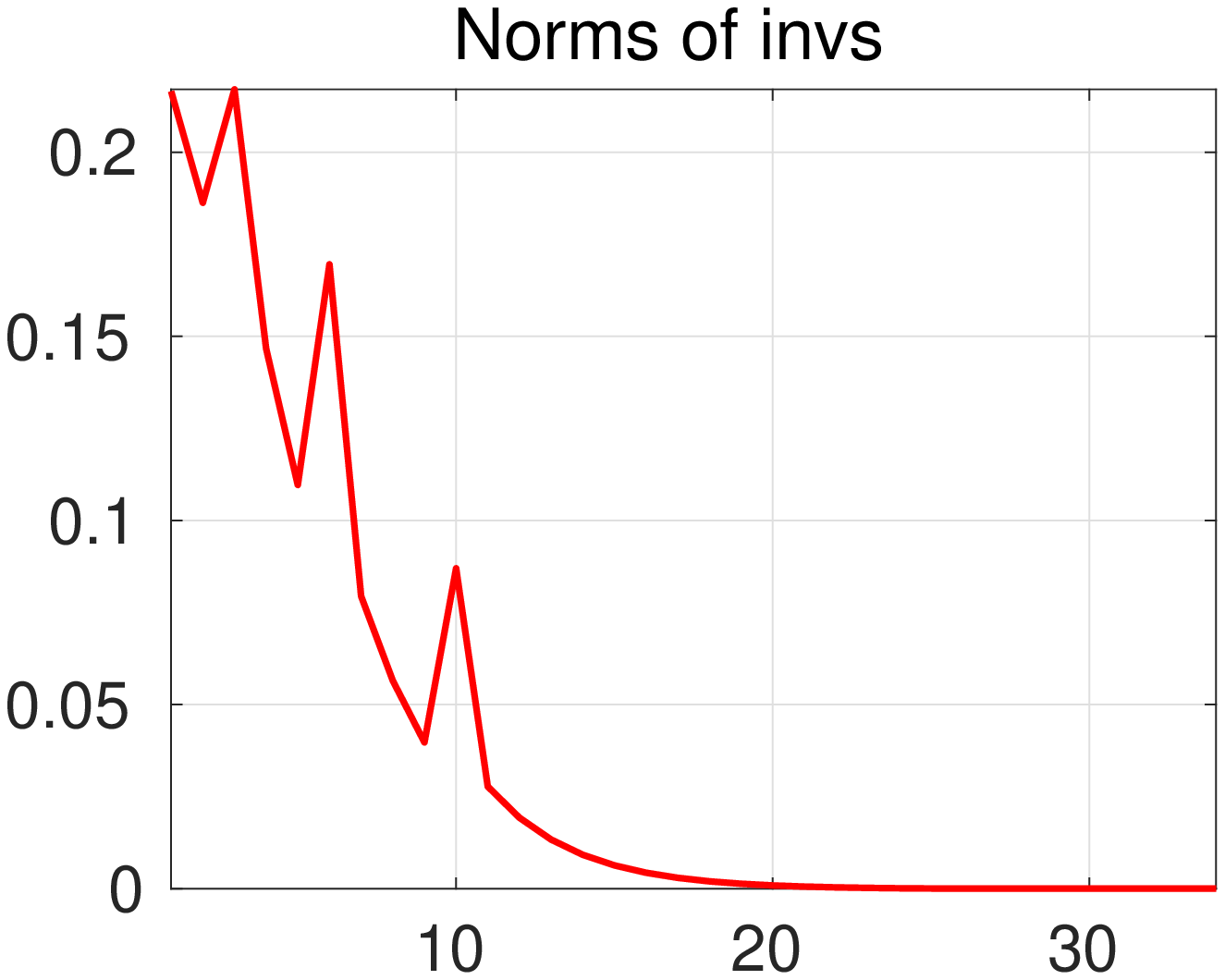}
\caption{\small Scaled norms  $\|{\bf Q}^{(1)}_k\|$ (left) and $\|{\bf Q}_0^{-1}\widehat{\bf Q}^{(1)}_k\|$
 (right) vs. $k=1,\ldots,R$.}
\label{fig:Scaled_Inv}
\end{figure}

We illustrate the ``add-and-compress" computational scheme in the following example.
 We consider the covariance matrix ${\bf C}_R$, obtained by 
 a rank-$R$ $\sinc$ approximation of the Slater function $e^{-\|x\|^p}$
 with $R=40$ on a grid with $n=1025$ sampling points. 
 Figure \ref{fig:Scaled_Inv} demonstrates  the decay in both the matrix 
 norms ${\bf Q}^{(1)}_k$ (left),  and the scaled, preconditioned matrices 
 ${\bf Q}_0^{-1}\widehat{\bf Q}^{(1)}_k$, $k=1,\ldots,R$ (right).
 We use the scaling factor of $1/n$.
 The right figure indicates that the analytic matrix functions ${\cal F}({\bf C}_R)$
 can be evaluated by using an exponentially fast convergent power series supported 
 by the ``add-and-compress'' strategy to control the tensor rank.
 
In the kriging calculations (Task 3 above), the low-rank tensor 
 structure in the covariance matrix ${\bf C}_R$
can be directly utilized if the sampling points in the kriging algorithm 
form a smaller $m_1\times m_2 \times m_3$ tensor sub-grid of the 
initial $n\times n \times n$ tensor grid $\Omega_n$ with $m_\ell < n$. 
The same argument also applies to the evaluation of conditional covariance.
In the general case of ``non-tensor'' locations of the sampling points,
some mixed tensor factorizations could be applied.

\subsection{Numerical illustrations}
\label{ssec:Numerics1}
In what follows, we check some examples of the low-rank Tucker tensor approximation of the $p$-Slater function
$C(x)=e^{-\|x\|^p}$, and Mat\'ern kernels
 with full-grid tensor representation. We demonstrate the fast exponential convergence 
 of the tensor approximation in the Tucker rank.
%
 The functions were sampled on the $n_1 \times n_2\times n_3$ 3D Cartesian grid  with
 $n_\ell=100$, $\ell=1,\,2, \,3$.
  \begin{figure}[htb]
\centering
\includegraphics[width=7.5cm]{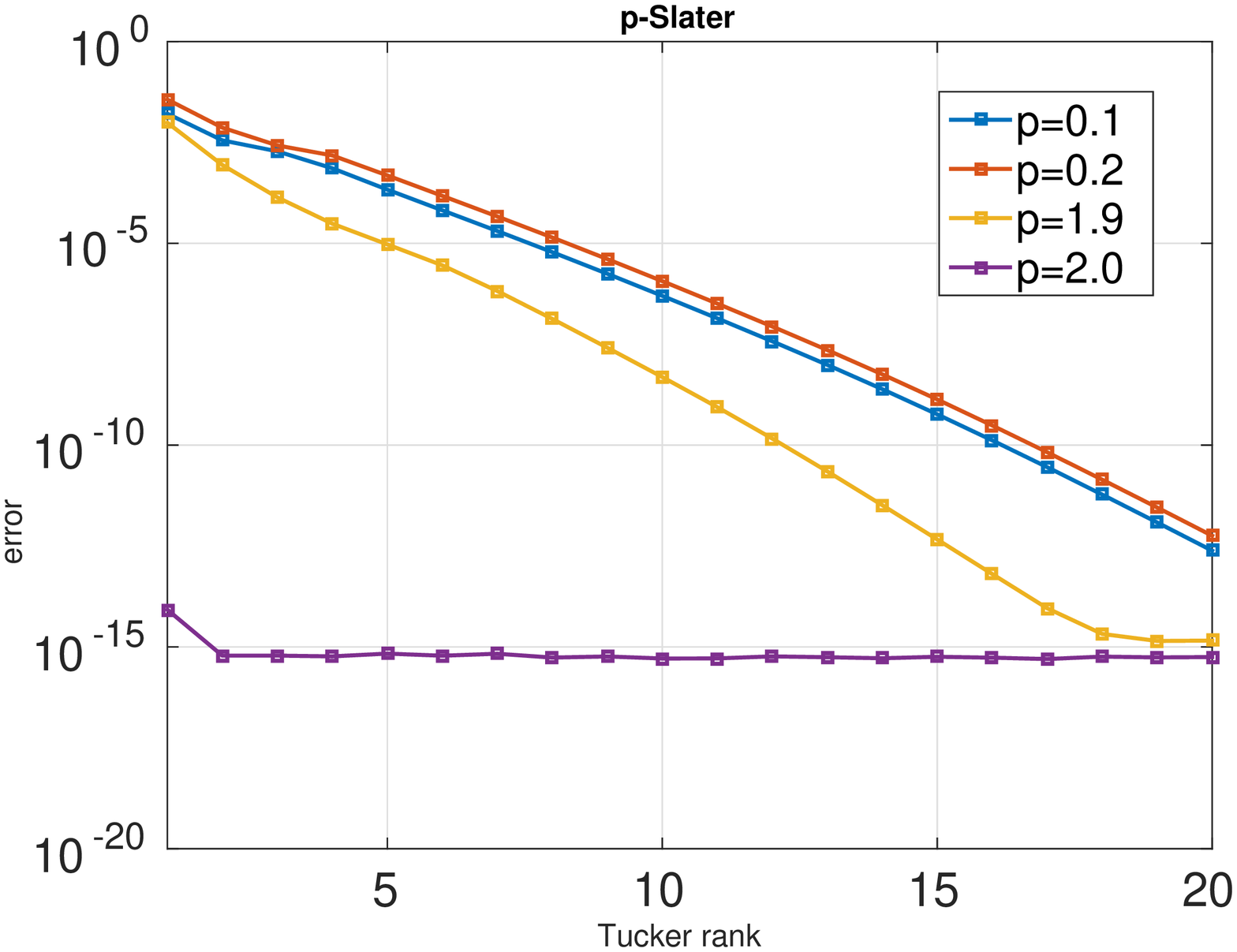}
\includegraphics[width=7.0cm]{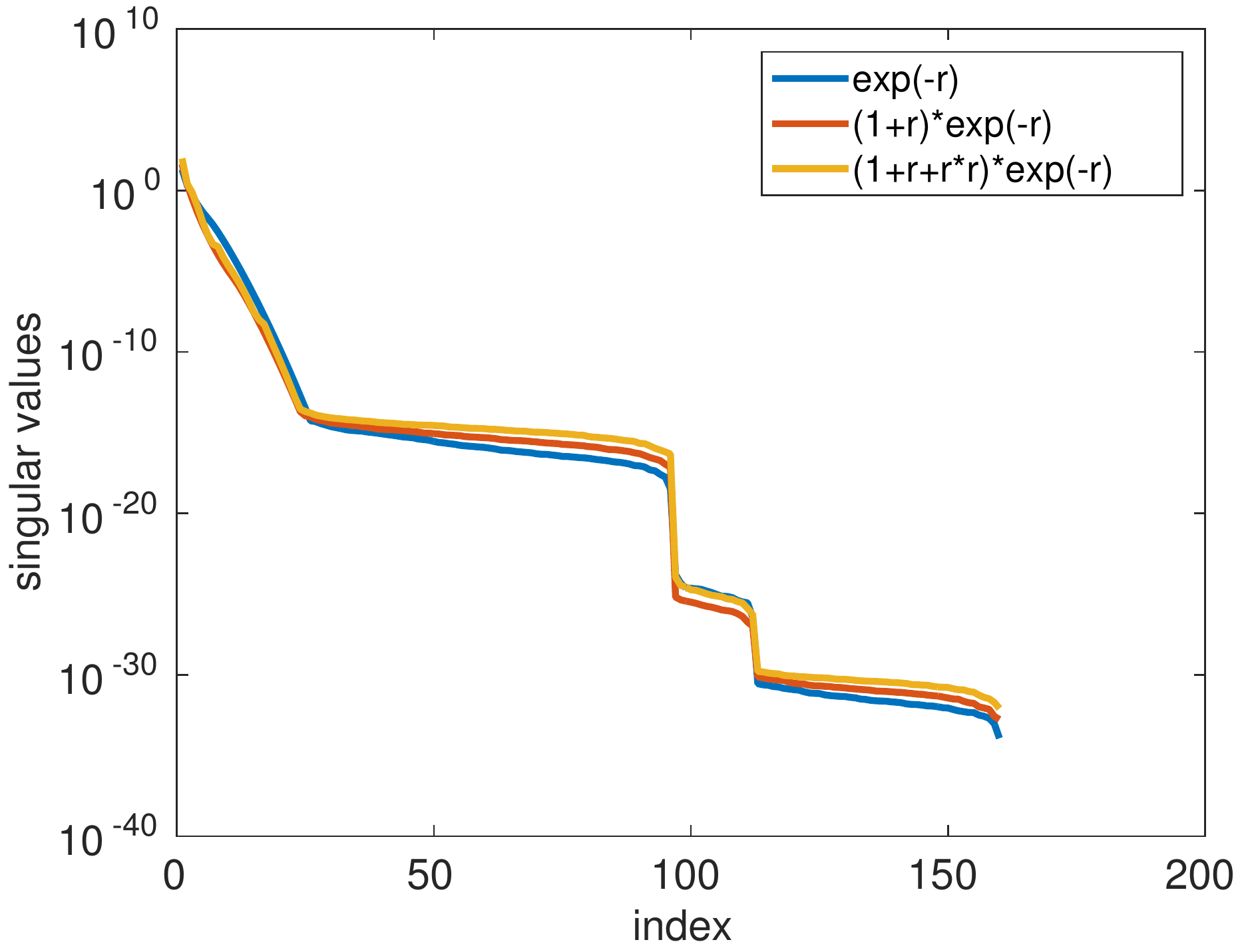}
\caption{\small Convergence in the Frobenius error (\ref{FNorm}) w.r.t. the Tucker rank for the
function (\ref{eq:exp_p}) with $p=0.1,\; 0.2,\; 1.9,\; 2.0$ (left); 
Decay of singular values of the weighted Slater function (right).}
\label{fig:Conv_Tucker_allp}
\end{figure}
  
For a continuous function $q:\Omega\to \mathbb{R}$, where
$\Omega := \prod^d_{\ell=1} [-b_\ell, b_\ell]\subset \mathbb{R}^d $, and
$0< b_\ell <\infty $, we introduce the 
collocation-type function-related tensor of order $d$:
$$
 {\bf Q}\equiv {\bf Q} (q) :=[q_{i_1 \ldots i_d}] \in
{\mathbb{R}}^{I_1 \times \ldots \times I_d} \mbox{ with }\;
q_{i_1 \ldots i_d} := q(x^{(1)}_{i_1},\ldots , x^{(d)}_{i_d}),
$$
where $(x^{(1)}_{i_1},\ldots ,x^{(d)}_{i_d})  \in \mathbb{R}^d$
are grid collocation points, indexed by ${\cal I} = I_1 \times \ldots \times I_d$, 
\begin{equation} \label{grid_points}
x^{(\ell)}_{i_\ell} = -b_\ell +(i_\ell -1) h_\ell, \quad
i_\ell = 1,2, \ldots , n_\ell,\; \ell=1,\ldots,d,
\end{equation}
which are the nodes of
equally spaced subintervals with a mesh size of $h_\ell=2b_\ell/(n_\ell-1)$.
%
%
  
 \begin{figure}[htb]
\centering
\includegraphics[width=7.6cm]{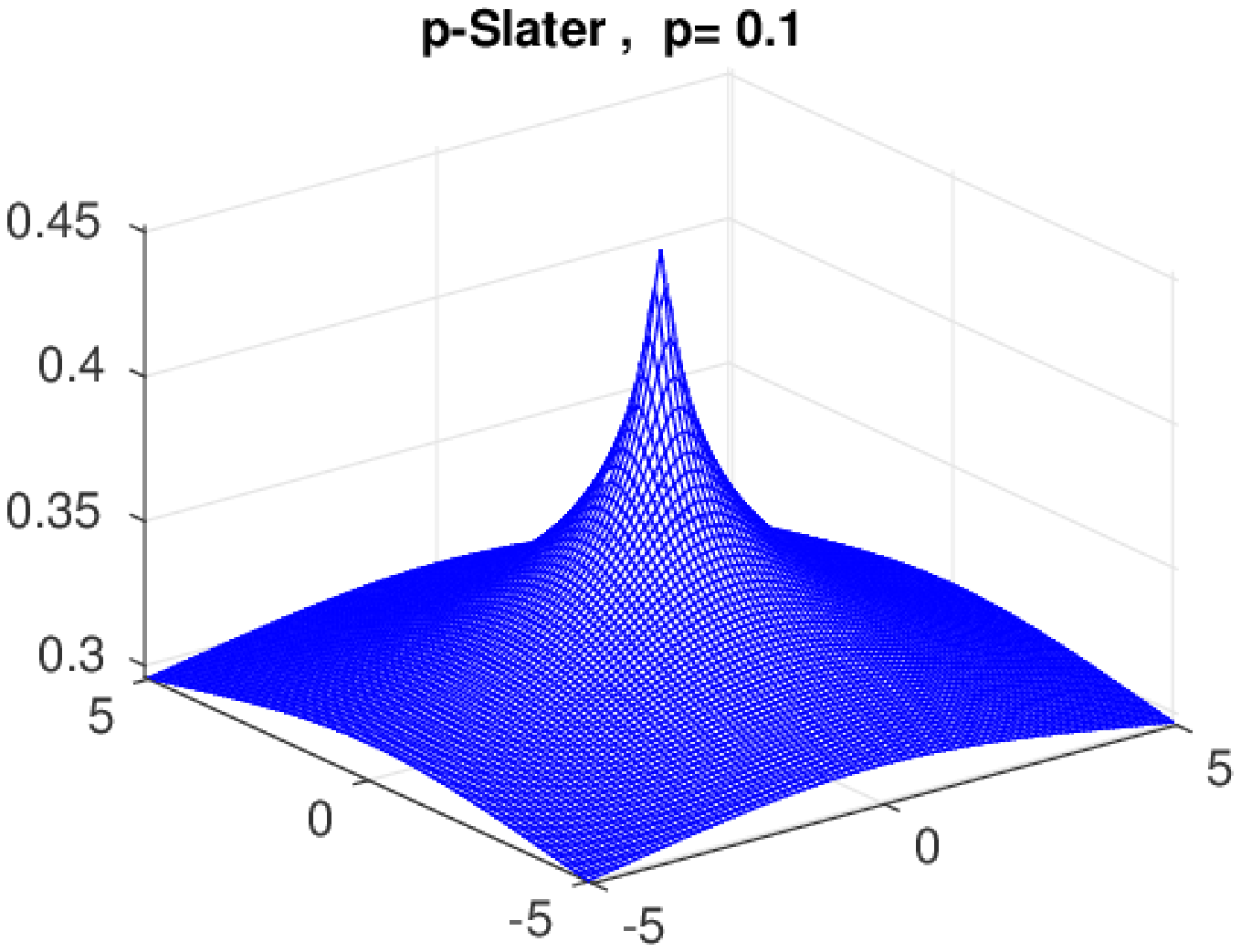} \quad
\includegraphics[width=7.6cm]{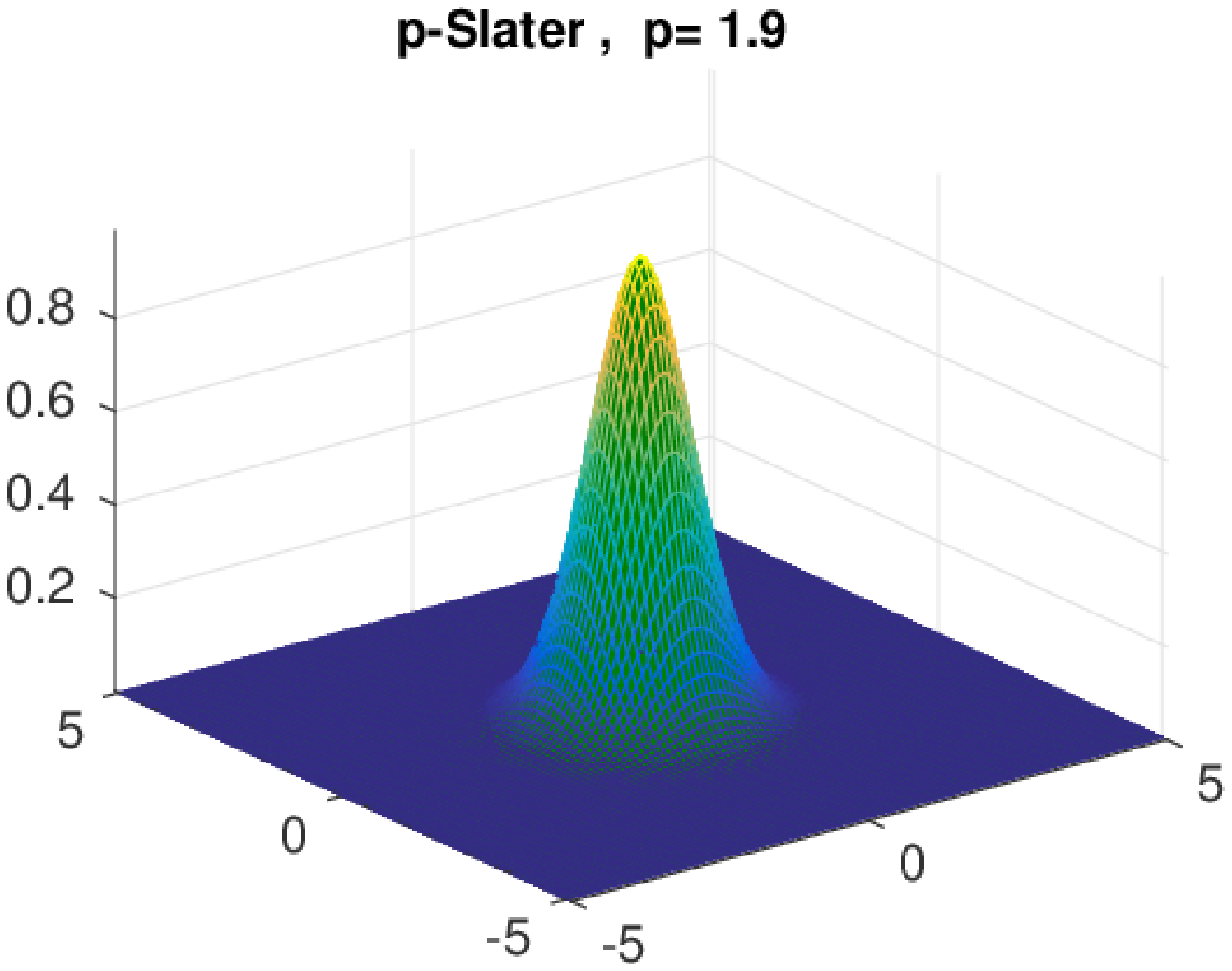}  
\caption{Cross section of the 3D radial function (\ref{eq:exp_p}) with $p=0.1$ (left) 
and $p=1.9$ (right) at level z=0.}
\label{fig:Surf_Slater_allp}
\end{figure}

 \begin{figure}[htb]
\centering
\includegraphics[width=8.6cm]{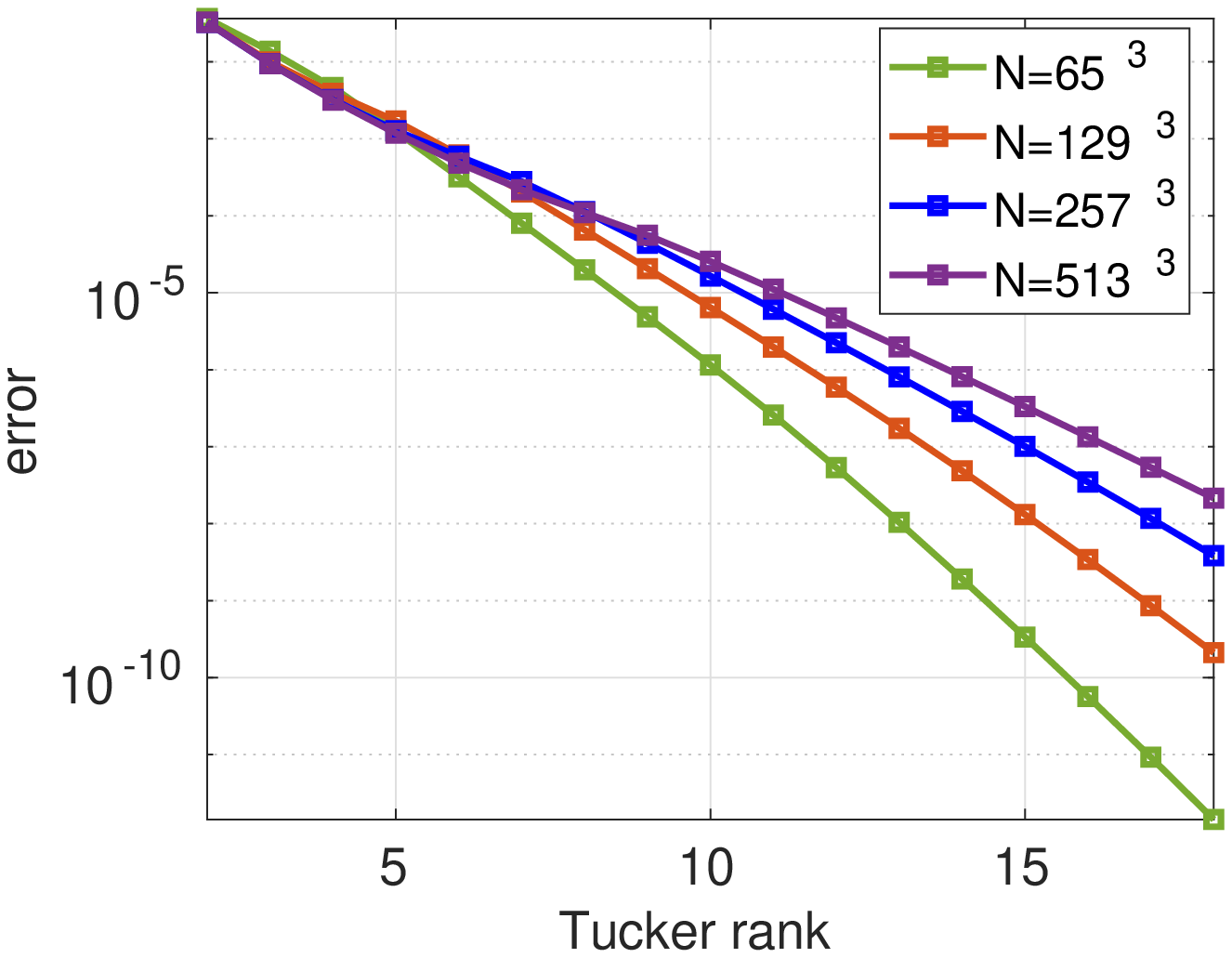}
\caption{Multigrid Tucker: convergence w.r.t. Tucker ranks of a Slater function 
with $p=1$ on a sequence of grids.}
\label{fig:Conv_Tucker_mg}
\end{figure}
We test the convergence of the 
error in the relative Frobenius norm with respect to the Tucker rank for $p$-Slater 
functions with $p= 0.1,\; 0.2, \; 1.9 ,\; 2.0$.  The Frobenius norm is computed as
 \begin{equation}
 \label{FNorm}
E_{FN} = \frac{\| {\bf Q} -{\bf Q}_{(r)}\|}
{||{\bf Q}||},
\end{equation}
 where ${\bf Q}_{(r)}$ is the tensor reconstructed from the Tucker rank-$r$ decomposition of
 ${\bf Q}$.
 
 Figure \ref{fig:Conv_Tucker_allp} shows convergence of the Frobenius error with respect to the
 Tucker  rank for the $p$-Slater function discretized on the 3D Cartesian grid 
 \begin{equation}\label{eq:exp_p}
 C(x,y)=e^{-\|x-y\|^p}
\end{equation}
 for different values of the parameter $p$. These 
 functions are illustrated in Figure~\ref{fig:Surf_Slater_allp} for 
 $p= 0.1$ and $p=1.9$.
 Figure \ref{fig:Conv_Tucker_mg}  
 shows for a Slater function with $p=1$ the dependence of the Tucker decomposition error
 in the Frobenius norm (\ref{FNorm}) on the Tucker rank, for the increasing grid parameter $n$.  
 
\begin{figure}[htb]
\centering
\includegraphics[width=7.2cm]{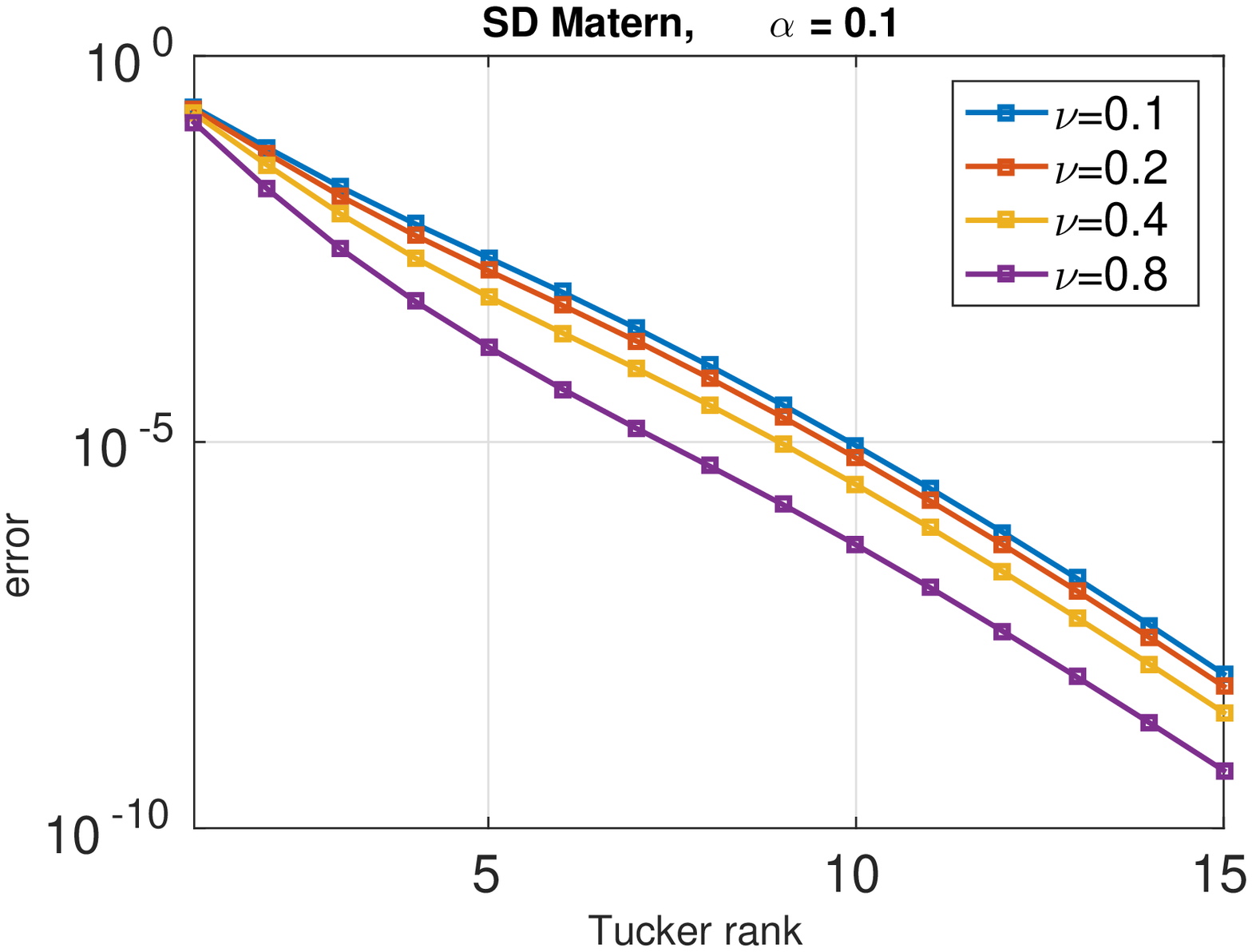} \quad
\includegraphics[width=7.0cm]{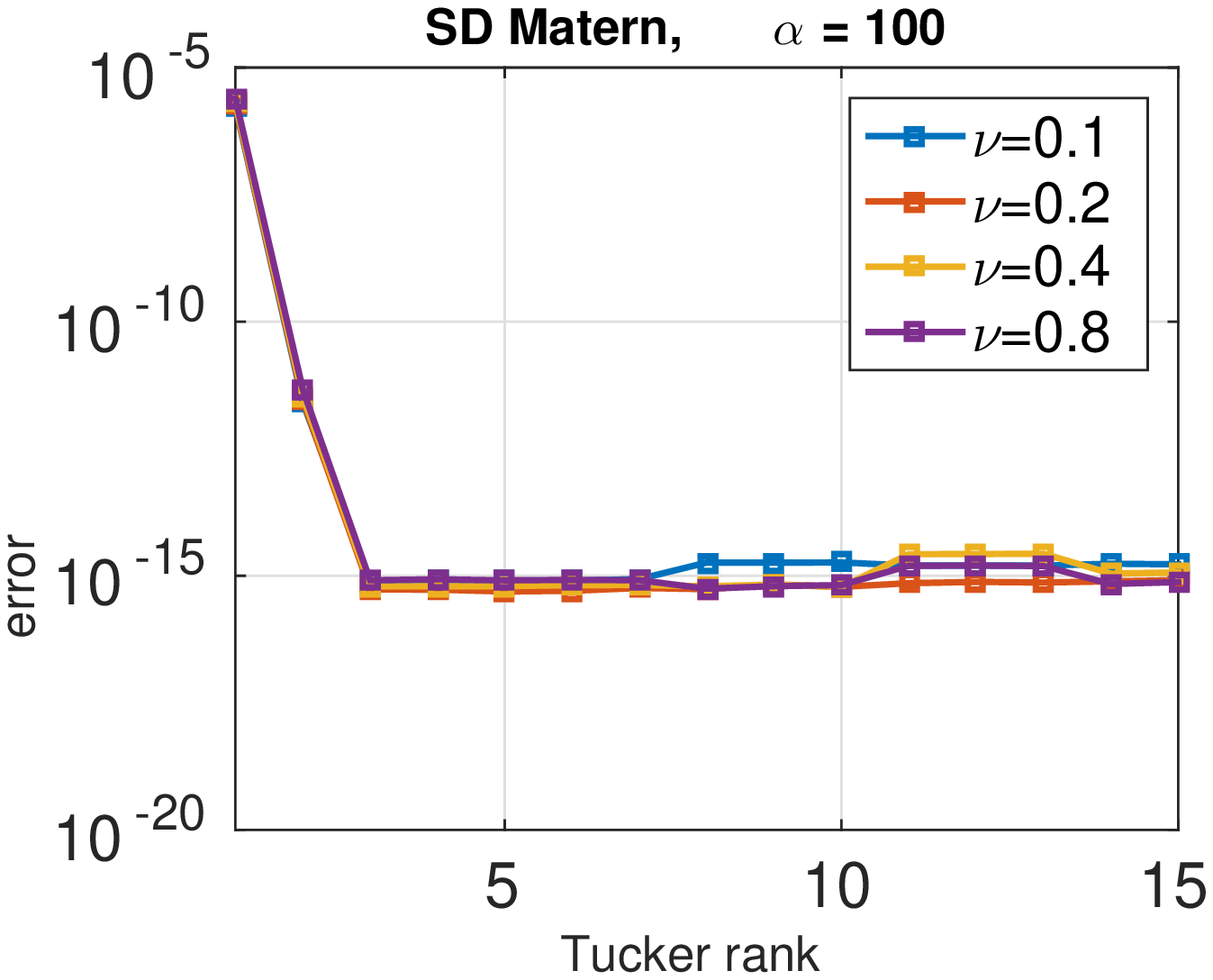}
\caption{Convergence w.r.t the Tucker rank of 3D spectral density of Mat\'ern covariance 
 (\ref{eqn:SD_matern}) with $\alpha=0.1$ (left) and $\alpha=100$ (right).}
\label{fig:Tuck_SD_Matern}
\end{figure}

\begin{figure}[htb]
\centering
\includegraphics[width=7.6cm]{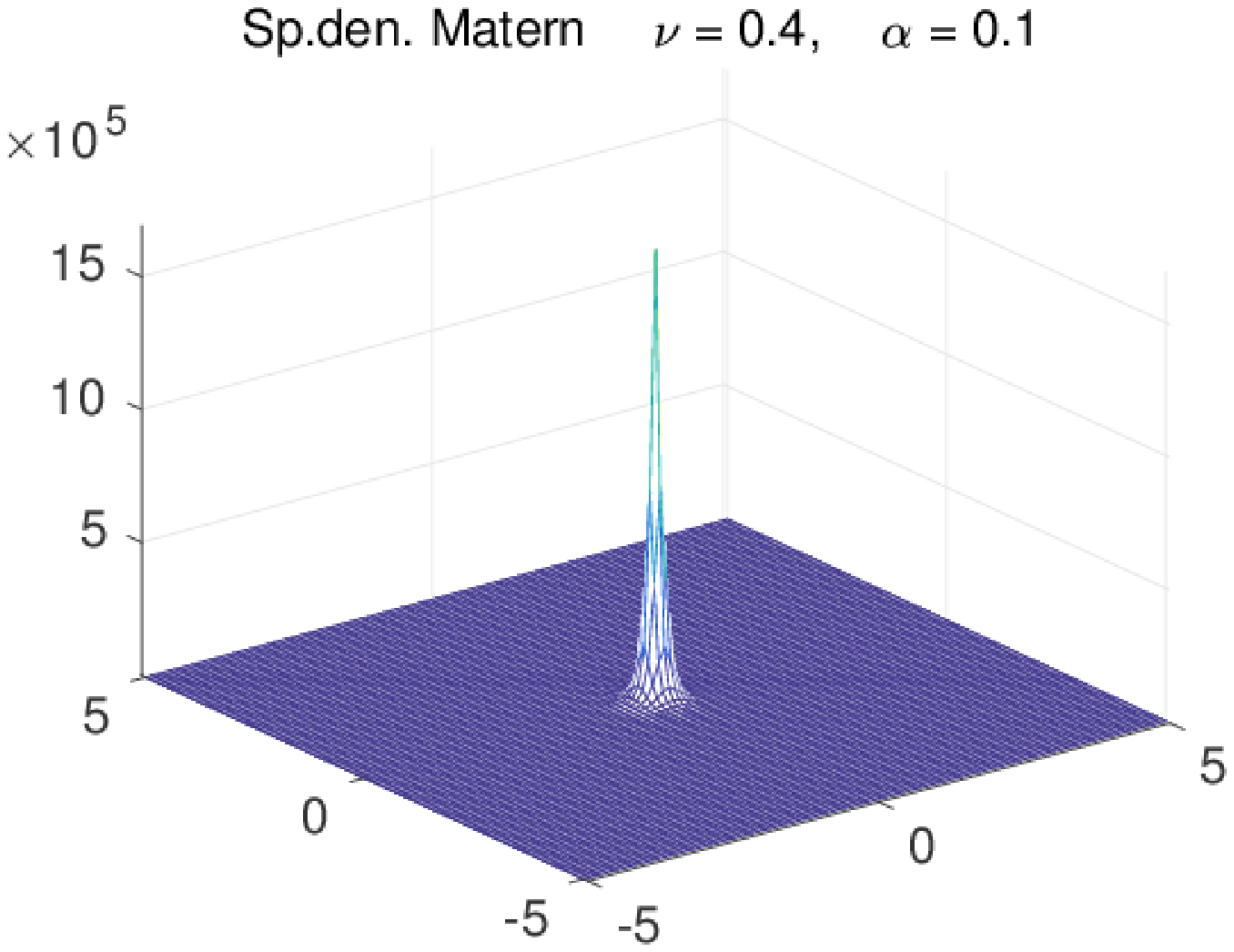} \quad
\includegraphics[width=7.6cm]{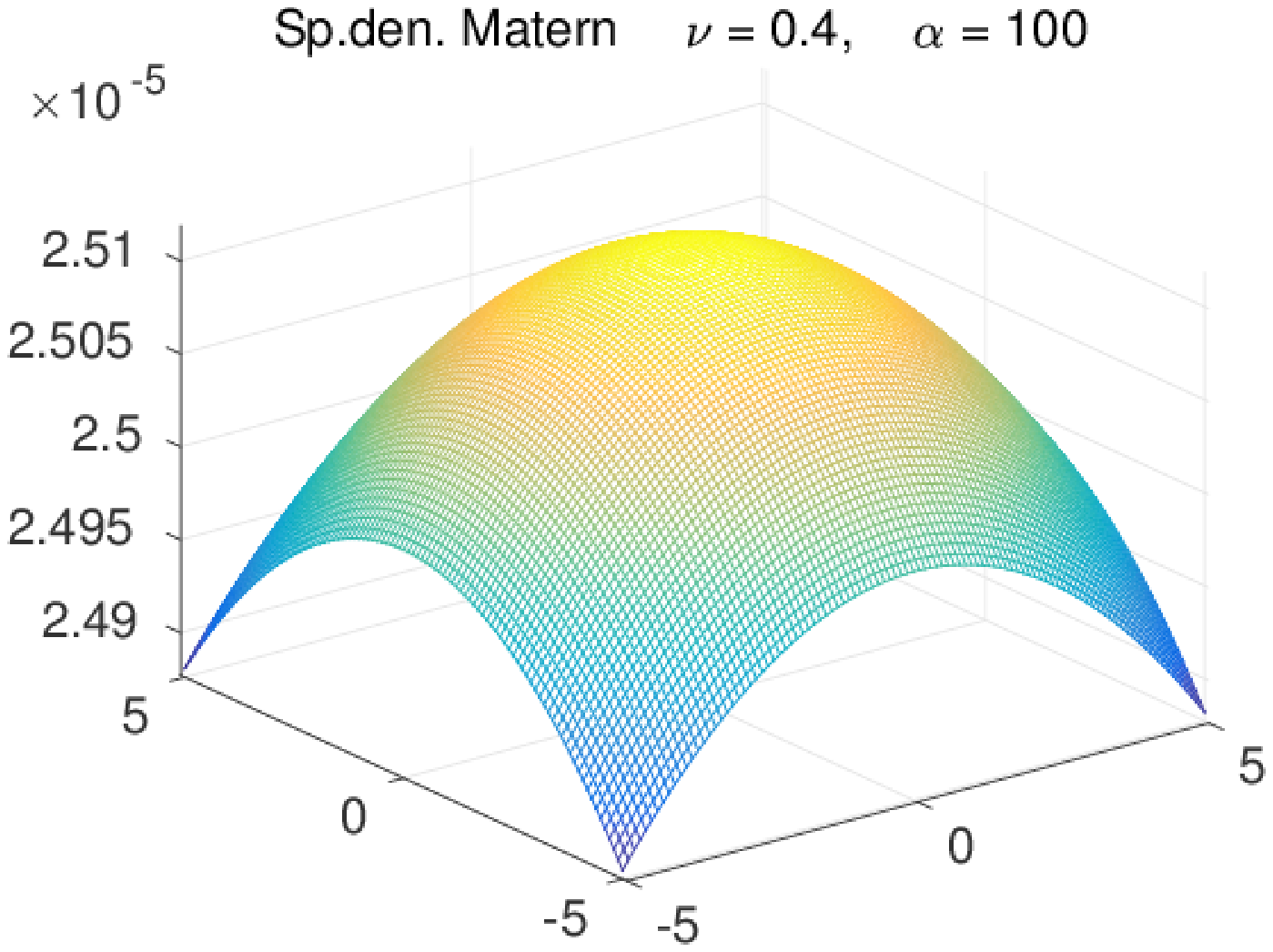}
\caption{The shape of 3D spectral density of Mat\'ern covariance 
 (\ref{eqn:SD_matern}) with $\alpha=0.1$ (left) and $\alpha=100$ (right).}
\label{fig:Surf_SD_Matern}
\end{figure}

 Figure \ref{fig:Tuck_SD_Matern} shows the convergence with respect to the Tucker rank
 for the spectral density of Mat\'ern covariance,
 \begin{equation}
\label{eqn:SD_matern}
f_{\alpha,\nu}(\rho):= \frac{C}{(\alpha^2+\rho^2)^{\nu+d/2}},
\end{equation}
where $\alpha \in (0.1,100)$ and $d=1,2,3$. The Tucker decomposition rank is strongly 
dependent on the parameter $\alpha$ and weakly depend on the
parameter $\nu$. The 3D Mat\'{e}rn functions
with the parameters $\nu=0.4$, $\alpha=0.1$ (left) and $\nu=0.4$, $\alpha=100$ (right) 
are presented in Fig.~\ref{fig:Surf_SD_Matern}, 
showing the function at $z=0$.

These numerical experiments demonstrate the good algebraic separability of the typical 
multidimensional functions used in spatial statistics, that lead us to apply the low-rank 
tensor decomposition methods to the multidimensional problems of statistical data analysis.

 
\section{Solutions to typical tasks in low-rank tensor format}
\label{sec:8tasksLow}

In this section, we walk through the solutions to the statistical
questions raised in the motivation above, Section~\ref{sec:Motivation}.
We add some lemmas to summarize the new computing and storage costs.
We let $N=n^d$, 
measurement vector $\bz \in \mathbb{R}^{m}$, $\bC_{ss}\in \mathbb{R}^{N\times N}$, 
$\bC_{zz}\in \mathbb{R}^{m\times m}$, and
$\bC_{sz}\in \mathbb{R}^{N\times m}$.
We also introduce the restriction operator $P$, which consists of only ones and zeros, 
and pick sub-indices $\{i_1,...,i_m\}$ 
from the whole index set $\{i_1,...,i_N\}$. This operator $P$ has tensor-1 structure, 
i.e., $P=\bigotimes_{\nu=1}^d P_{\nu}$. 
An application of this restriction tensor does not change tensor ranks.\\
\textbf{Computing matrix-vector product:}
We let $\bC= \sum_{i=1}^{r} \bigotimes_{\mu=1}^d \bC_{i \mu}$. If $\bz$ is separable, 
i.e., $\Vert \bz - \sum_{j=1}^{r_b} \bigotimes_{\nu=1}^d \bz_{j \nu}\Vert \leq \varepsilon $, then 
\begin{equation}
 \bC \bz= \sum_{i=1}^{r}\sum_{j=1}^{r_z}\bigotimes_{\mu=1}^d \bC_{i \mu} \bz_{j \mu}.
\end{equation}
If $\bz$ is non-separable, then low-rank tensor properties cannot be employed and either 
the FFT idea \cite{nowak2013kriging} or the hierarchical matrix technique should be 
applied instead \cite{khoromskij2009application,HackHMEng,Part1,Part2}.
\begin{lemma}
\label{lem:MVprod}
The computing cost of the product $ \bC \bz$ is reduced from $\mathcal{O}(N^2)$ 
to $\mathcal{O}(rr_zdn^2)$, where $N=n^d$, $d\geq 1$.
\end{lemma}
\textbf{Trace and diagonal of $\bC$:}
Let $\bC \approx \tilde{\bC} =  \sum_{i=1}^r \bigotimes_{\mu=1}^d \bC_{i \mu}$, then
\begin{equation}
\label{eq:diagC}
\diag(\tilde{\bC}) = \diag\left( \sum_{i=1}^r \bigotimes_{\mu=1}^d \bC_{i \mu}\right )=
\sum_{i=1}^r \bigotimes_{\mu=1}^d \diag\left( \bC_{i \mu}\right ),
\end{equation}
\begin{equation}
\label{eq:traceC}
\trace(\tilde{\bC}) = \trace\left(\sum_{i=1}^r \bigotimes_{\mu=1}^d \bC_{i \mu} \right)=
\sum_{i=1}^r \prod_{\mu=1}^d \trace(\bC_{i \mu}). 
\end{equation}
The proof follows from the properties of the Kronecker tensors.
\begin{lemma}
The cost of computing the right-hand sides in Eq.~\ref{eq:diagC}-Eq.~\ref{eq:traceC} 
is $rdn$, where $\bC_{i \mu} \in\mathbb{R}^{n\times n}$. 
\end{lemma}
For simplicity, 
we assume that $n_1=n_2=...=n_d=n$ and $\sum_{i=1}^d n_i=dn$.
\begin{lemma}
\label{lem:tracediag}
The computing cost of $\diag(\bC)$ and $\trace(\bC)$ is reduced from $\mathcal{O}(N)$ to 
$\mathcal{O}(rdn)$. The cost of $\det(\bC)$ is reduced from $\mathcal{O}(N^3)$ to 
$\mathcal{O}(dn^3)$.
\end{lemma}
\begin{ex}
A simple Matlab test for computing $\trace(\bC)$ on a working station with 128GB produces 
the computing times in seconds shown in Table~\ref{table:trace_comp_time}.
 \begin{table}[htbp]
 \centering
 \begin{tabular}{|c|c|c|c|}
 \hline
                  &  $n=100$  &  $n=500$   & $n=1000$   \\  \hline
   $d=1000$   &   3.7       & 67& 491\\ \hline
 \end{tabular}
 \caption{Computing time (in sec.) to set up and compute the trace of $ \tilde{\bC} = 
 \sum_{j=1}^{r} \bigotimes_{\nu=1}^d \mathbf{\bC}_{j \nu} $, 
 $r=10$, $ \tilde{\bC} \in \mathbb{R}^{N\times N}$, 
 where $N=n^d$, $d=1000$ and  $n=\{100, 500, 1000\}$. A modern desktop computer with 128 GB RAM was used.}
 \label{table:trace_comp_time}  
\end{table}
\end{ex}

In what follows, we discuss the computation of trace and diag in the case of Tucker representation of the 
generating Mat{\'e}rn function on $(2k+1)^{\otimes d}$ grid and present the corresponding numerical example
for $d=3$. We notice that due to Toeplitz structure of the skeleton matrices of size 
$(2k+1)\times (2k+1)$, the diagonal of the covariance matrix ${\bf C} $ is the weighted Tucker 
sum of the Kronecker products of scaled identity matrices in $\mathbb{R}^{2k+1}$. The scaling factor 
is determined by the value of generating Tucker vector ${\bf v}_{\nu_\ell}^{(\ell)}(k+1)$
corresponding the origin of the computational box $[-b,b]^d$. 
Let the matrix ${\bf C} $ be composed by using the Tucker tensor ${\bf A} $ approximating the 
the generating Mat{\'e}rn function.
Then, as a straightforward consequence of the above remark, we derive the simple representations
\[
\mbox{diag}({\bf C})=  {\bf A}(x_0) {\bf E}^{(d)}, \quad 
 \mbox{trace}({\bf C})= {\bf A}(x_0) (2k+1)^d,
\]
where $x_0$ corresponds to the origin $x=0$ in $[-b,b]^d$ and ${\bf E}^{(d)}$ is the identity 
matrix in the full tensor space. The grid coordinate of $x_0$ is determined by the multi-index
$(k+1,...,k+1)$.

Next table represents the values of ${\bf A}(x_0)$ computed by the Tucker approximation to
the 3D Slater function  $e^{-\|x \|}$ on the $(2k+1)^{\otimes 3}$ grids with $n=2k+1=129, 257, 513$, 
and for different Tucker rank parameters $r=1,2,...,10$. 

\begin{table}[htbp]
 \begin{center}
 \begin{tabular}
[c]{|c|c|c|c|c|c|c|c|c|c|c|c| }%
\hline
Tucker rank $r$ & 1   &   2  &   3    & 4   & 5   & 6  &  7 &  8 &  9 &  10 \\
\hline
n=129   & 0.386   & 0.20 &  0.12 & 0.07 & 0.04 & 0.017 & 0.002 &  1.2e-4 & 8.4e-6 & 7.5e-6 \\
\hline
n=257   & 0.386   & 0.20 &  0.12 & 0.073 & 0.046 & 0.029 & 0.017 &  0.007 & 8.0e-4 & 1.4e-5 \\
\hline
n=513   & 0.386   & 0.20 &  0.12 & 0.073 & 0.047 & 0.031 & 0.020 &  0.0138 & 0.008 & 0.0035 \\
 \hline
 \end{tabular}
\end{center}
\caption{The error of the Tucker approximation to the value  ${\bf A}$ at the origin
(the exact value is equal to $1$) versus the Tucker rank and the grid size $n=2k+1$.}
\label{tab:zero_points}
\end{table}

Given the rank-${\bf r}$ Tucker tensor  ${\bf A} $, the complexity for calculation of 
${\bf A}(x_0)$ is estimated by $O(r^d)$.

\textbf{Computing square root $\bC^{1/2}$:}
$\bC^{1/2}$ can be computed as in \eqref{eq:sqrt}. An iterative method for computing $\bC^{1/2}$ 
is presented in \cite{gavrilyuk2005data}.\\
 
\textbf{Linear solvers in a low-rank tensor format:} Likely, there is already a good theory for solving  
linear systems $ \bC\bw= \bz$ 
with symmetric and positive definite matrix $\bC$, in a tensor format. We refer to the overview works
\cite{KhorTensorReview,khor-survey-2014, GraKresTo:13,hackbusch2012tensor}. Some particular linear 
solvers are developed in 
\cite{DoKhSa-qtt_fft-2011,BallaniKressner,Khoromskij_Low_Kronecker,Khoromskij_Low_Tacker,DolgLitv15}, 
\cite{dolgov2014computation, Khoromskij_Low_Tacker}. 
We also recommend to use the Tensor Toolbox \cite{Kolda:07}, which contains routines for CP 
and Tucker tensor formats. 
 
\subsection{Computing $\bz^T \bC^{-1} \bz$ }
\begin{lemma}
\label{lem:lowrank_quad}
We let $\Vert \bz - \sum_{i=1}^{r} \bigotimes_{\mu=1}^d {\bz}_{i \mu}\Vert \leq \varepsilon$. 
We assume that there is an iterative method 
that can be used to solve the linear system $ \bC\bw= \bz$ in a low-rank tensor format and to 
find the solution in the  
form $\bw = \sum_{i=1}^{r} \bigotimes_{\mu=1}^d {\bw}_{i \mu}$.
Then, the quadratic 
form $\bz^T \bC^{-1} \bz$ is the following scalar products:
\begin{equation}
\bz^T \bC^{-1} \bz= \sum_{i=1}^{r}\sum_{j=1}^{r_z} \prod_{\mu=1}^d   
\left ( \bw_{i \mu}, \bz_{j \mu} \right),
\end{equation}
\end{lemma}
The proof follows from the definition and properties of the tensor and scalar products.\\
If $\bz$ is non-separable, then low-rank tensor properties cannot be employed and the 
FFT idea \cite{nowak2013kriging} or the hierarchical 
matrix technique \cite{khoromskij2009application} should be employed.
\begin{lemma}
\label{lem:quadratic}
The computing cost of the quadratic form $\bz^T \bC^{-1} \bz$ is the product of the number of required iterations
and the cost of one iteration, which is
$\mathcal{O}(rr_zd m^2)$ (assuming that the iterative method required only matrix-vector products).
\end{lemma}
The proof follows from the definitions and properties of the tensor and scalar products.
\subsection{Interpolation by simple kriging}
The three most computationally demanding
tasks in kriging are: (1) solving an $M\times M$ system of equations to obtain the
kriging weights; (2) obtaining the $N\times 1$ kriging estimate by
superposing the kriging weights with the $N\times M$ cross-covariance
matrix between the measurements and the unknowns; and (3) evaluating the $N \times 1$ estimation
variance as the diagonal of an $N\times N$ conditional covariance matrix \cite{nowak2013kriging}. 
Here, $M$ refers
to the number of measured data values, and $N$ refers to the number of estimation
points.
When optimizing the design of the sampling patterns,
the challenge is to evaluate the scalar measures of the $N \times N$ conditional 
covariance matrix (see $\phi_A$ and $\phi_C$ 
in Eq.~\ref{eq:Ccriterion} and Task 4)
repeatedly within a high-dimensional and non-linear optimization procedure
(e.g., \cite{Kollat_al_2008AWR_epsilon_hBOA_MOO, Shah_Reed_2011OR_MOO_comparison_knapsack}).

The following kriging formula is well-known \cite{nowak2013kriging}:
\begin{equation}
\hat{\mathbf{s}}=\mathbf{C}_{sz}\mathbf{C}_{zz}^{-1}\mathbf{z}.
\label{eq:estimate}
\end{equation}
\begin{lemma}
If
$\Vert \bC_{sz} - \sum_{i=1}^{r_C} \bigotimes_{\mu=1}^d {\bC}_{i \mu}\Vert \leq \varepsilon$, 
for some small $\varepsilon\geq 0$, and
Lemma~\ref{lem:lowrank_quad} holds, then 
\begin{equation}
\label{eq:lrkri}
\mathbf{C}_{sz}\mathbf{C}_{zz}^{-1}\mathbf{z}\approx 
\sum_{i=1}^{r_z}\sum_{j=1}^{r_C}\bigotimes_{\nu=1}^d  
{\bC}_{i \nu} {\bw}_{j \nu}.
\end{equation}
\end{lemma}
The proof follows from the definitions and the properties of the tensor and scalar products.\\
\begin{lemma}
\label{lem:quadratic}
The computing cost of solving the linear system $\mathbf{C}_{zz}^{-1}\mathbf{z}$ 
is $\mathcal{O}(\#\mbox{iters}\cdot r_zrd m^2)$. 
Computation of the kriging coefficients by Eq.~\ref{eq:lrkri} costs 
$\mathcal{O}(r_zr_Cd nm)+\mathcal{O}(\#\mbox{iters}\cdot r_zrd m^2)$.
\end{lemma}
If $\bz$ is non-separable, then the low-rank tensor properties cannot be employed and 
either the FFT idea \cite{nowak2013kriging} or the hierarchical matrix 
technique \cite{khoromskij2009application} should be applied.
%
%
 %
\subsection{Computing conditional covariance}
We let $\by \in \mathbb{R}^{m}$ be the vector of measurements.
The conditional covariance matrix is
\begin{equation}
\mathbf{C}_{ss \vert y}
= \mathbf{C}_{ss}           - \mathbf{C}_{sy} \mathbf{C}_{yy}^{-1}\mathbf{C}_{ys}.
\end{equation}
The associated estimation variance $\hat{\bm{\sigma}}$ is
the diagonal of the $N \times N$ conditional covariance matrix $\mathbf{C}_{ss|y}$: 
\begin{align}
\hat{\boldsymbol{\sigma}}_{\mathbf{s}} &= \diag(\mathbf{C}_{ss|y})
=\diag\left( \mathbf{C}_{ss}           - \mathbf{C}_{sy} \mathbf{C}_{yy}^{-1}\mathbf{C}_{ys} 
\right).
\label{eq:estvar}
\end{align}
We assume that 
the measurements are taken at locations that form a subset of the total 
set of nodes $\mathcal{I}=\{0,...,N-1\}$;
i.e., $\mathcal{I}_{\mathcal{M}}=\{i_{1},...,i_{m}\}\subset \mathcal{I}$. 
We also assume that the nodes $\mathcal{I}_{\mathcal{M}}$ belong 
to a tensor mesh, i.e., if $\mathcal{I}=\bigotimes_{\nu=1}^d I_{\nu}$, 
and $\mathcal{I}_{\mathcal{M}}=\bigotimes_{\nu=1}^d \hat{I}_{\nu}$, then $\hat{I}_{\nu} \subseteq I_{nu}$.

We let
$\mathbf{C}_{yy} = \sum_{k=1}^{r} \bigotimes_{\mu=1}^d \bC_{k \mu}$. We again use 
low-rank tensor solvers, this time to solve the matrix system
$\mathbf{C}_{yy}\bW = \mathbf{C}_{ys}$. We obtain the solution 
$\bW= \mathbf{C}_{yy}^{-1} \mathbf{C}_{ys}=\sum_{j=1}^{r_w} \bigotimes_{\mu=1}^d \bC_{j \mu}$.
Then, assuming that $\bC_{sy} \approx \sum_{i=1}^{r_C} \bigotimes_{\mu=1}^d {\bC}_{i \mu}$, we obtain
\begin{equation}
\mathbf{C}_{sy}\bW = \mathbf{C}_{sy}\mathbf{C}_{yy}^{-1} \mathbf{C}_{ys} = 
\sum_{i=1}^{r_C} \bigotimes_{\nu=1}^d \bC_{i \nu}
\sum_{j=1}^{r_w} \bigotimes_{\mu=1}^d \bC_{j \mu}=
\sum_{i=1}^{r_C}\sum_{j=1}^{r_w} \bigotimes_{\mu=1}^d \bC_{i \mu}\bC_{j \mu}
\approx \sum_{j=1}^{r_0}
\bigotimes_{\mu=1}^d \tilde{\bC}_{j \mu}, 
\end{equation}
where $r_0>0$ is the new rank after a rank-truncation procedure and $ \tilde{\bC}_{j \mu}$ are new factors.
We let $\mathbf{C}_{ss} = \sum_{i=1}^{r_s} \bigotimes_{\mu=1}^d \bC_{i \mu}$.
The conditional covariance is
\begin{equation}
\mathbf{C}_{ss \vert y}=\sum_{i=1}^{r_s} \bigotimes_{\mu=1}^d \bC_{i \mu} - 
\sum_{i=1}^{r_0} \bigotimes_{\mu=1}^d \tilde{\bC}_{i \mu} =\sum_{i=1}^{r_s+r_0} 
\bigotimes_{\mu=1}^d \hat{\bC}_{i \mu}, 
\end{equation}
where $\hat{\bC}_{i \mu}=\bC_{i \mu}$ for $1\leq i\leq r_s$ and 
$\hat{\bC}_{i \mu}=-\tilde{\bC}_{i \nu}$ for $r_s < i \leq r_0+r_s$.
%

%
 \subsection{Example: Separable covariance matrices}
\label{sec:ExampleSeparable}
We let $\cov(\bx,\by)=\exp^{-\vert \bx-\by\vert^2}$ be the Gaussian covariance function, 
where $\bx=(x_1,..,x_d)$, and $\by=(y_1,...,y_d)\in \mathcal{D} \subset \mathbb{R}^d$.
$\cov(\bx,\by)$ can be written as a tensor product 
of 1D functions: 
$$\cov(\bx,\by)=\exp^{-\vert x_1-y_1\vert^2}\otimes \ldots
\otimes \exp^{-\vert x_d-y_d\vert^2}.$$
After discretization of $\cov(\bx,\by)$, we obtain $\bC$ as a rank-1 Kronecker product 
of the 1D covariance matrices, i.e.,  
\begin{equation}
\label{eq:Crank1}
\bC=\bC_1 \otimes ... \otimes \bC_d.
\end{equation}
We note that arbitrary discretization (anisotropy) can occur in any direction.
\begin{lemma}
\label{thm:d-Chol-Kronecker}
If $d$ Cholesky decompositions exist, i.e, $\bC_i=\bL_i\cdot \bL_i^T$, and $i=1,...,d$, then 
\begin{equation}
\bC_1 \otimes ... \otimes \bC_d=(\bL_1\bL_1^T)\otimes ... \otimes (\bL_d \bL_d^T)= 
(\bL_1 \otimes ... \otimes \bL_d)\cdot (\bL_1^T \otimes ... \otimes \bL_d^T)=:\bL\cdot \bL^T,
\end{equation}
where $\bL:=\bL_1 \otimes...\otimes \bL_d$ and $\bL^T:=\bL_1^T \otimes...\otimes \bL_d^T$ 
are also lower- and upper-triangular matrices, respectively.
 \end{lemma}
Lemma \ref{thm:d-Chol-Kronecker} shows that a) the Gaussian covariance function in 
dimensions $d>1$ can be written as the tensor sum of 1D covariance 
functions, and b) its Cholesky factor can be computed via Cholesky factors 
computed from 1D covariances. The computational complexity drops 
from $\mathcal{O}(N\log N)$, $N=n^d$, to  $\mathcal{O}(d n\log n)$, where $n$ 
is the number of mesh points in 1D problem. Further research is required 
on non-Gaussian covariance functions.
\begin{lemma}
\label{thm:d-Kronecker-inv}
We let $\bC=\bC_1 \otimes ... \otimes \bC_d$.
If the inverse matrices $\bC_i^{-1}$, $i=1..d$, exist, then 
\begin{equation}
(\bC_1 \otimes ... \otimes \bC_d)^{-1} = \bC_1^{-1}\otimes ... \otimes \bC_d^{-1}.
\end{equation}
 \end{lemma}
The computational complexity drops from $\mathcal{O}(N\log N)$, $N=n^d$, to  
$\mathcal{O}(d n\log n)$, where $n$ is the number of mesh points in a 1D problem. 
\begin{remark}
\label{rem:compl}
We assume here that we have an efficient method to invert $\bC$ (e.g., FFT or 
hierarchical matrices) with a cost of $\mathcal{O}(N\log N)$. 
If not, then the complexity cost drops from $\mathcal{O}(N^3)$ to  
$\mathcal{O}(d n^3)$ (usual Gaussian elimination algorithm).
\end{remark}

\begin{lemma} 
\label{thm:d-det}
If $\bC_i$, $i=1..d$, are covariance matrices, then we can compute $\log\det \bC$, 
where $\bC$ is a separable rank-1 $d$-dimensional covariance function, as
\begin{equation}
\label{eq:tensor_logdet}
\log \det(\bC_1\otimes \ldots \otimes \bC_d)=\sum_{j=1}^d \log \det \bC_j \prod_{i=1, i\neq j}^d n_i.
\end{equation}
\end{lemma}
Proof: We check for $d=2$: $\det(\bC_1\otimes \bC_2)=\det(\bC_1)^{n_2} \cdot \det(\bC_2)^{n_1}$,
and then apply mathematical induction.

The computational cost 
drops again from $\mathcal{O}(N\log N)$, $N=n^d$, to $\mathcal{O}(dn\log n)$. 
A similar assumption to Remark~\ref{rem:compl} 
for computing $\det(\bC)$ also holds here.\\
\begin{ex}
We let $n=6000$, $d=3$, and $N=6000^3$. 
Using MATLAB on a MacBookPro with 16 GB RAM, the time required set up the matrices 
$\bC_1$,$\bC_2$, and $\bC_3$ is 11 seconds; it takes 4 seconds to compute $\bL_1$, 
$\bL_2$, and $\bL_3$. The large matrices $\bC$ and $\bL$ 
are never constructed (i.e., the Kronecker product is never calculated).
\end{ex}
\begin{ex}
In previous work \cite{LitvHLIBpro} we used the hierarchical matrix technique to 
approximate $\bC_{i}$ and its Cholesky 
factor $\bL_{i}$ for $n=2\cdot 10^6$ in 2 minutes. 

Here, we combine the hierarchical matrix technique and the Kronecker tensor product.

Assuming $\bC=\bC_1\otimes \ldots \otimes \bC_d$, we approximate $\bC$ for 
$n=(2\cdot 10^6)^d$ in $2d$ minutes. 
\end{ex}

\begin{lemma}

We let $\bC$ be the same as in Eq.~\ref{eq:Crank1}, and we let Lemma~\ref{lem:lowrank_quad} hold.
Then, we apply the property in Eq.~\ref{eq:tensor_logdet}
to obtain a tensor approximation of the log-likelihood:
 \begin{equation}
\label{eq:tensor_Log}
\LL \approx\tilde{\LL} = -\frac{\prod_{\nu=1}^d n_{\nu}}{ \log(2\pi)} - \sum_{j=1}^d 
\log \det \bC_j \prod_{i=1, i\neq j}^d n_i -\sum_{i=1}^r \sum_{j=1}^r 
\prod_{\nu=1}^d (\bu_{i,\nu}^T, \bu_{j,\nu}).
\end{equation}
\end{lemma}
Equation~\ref{eq:tensor_Log} shows one disadvantage of the Gaussian log-likelihood function in high dimensions. 
Namely, the log-likelihood grows exponentially with $d$ as $n^d$.

\section{Conclusion}

In this work, we demonstrate that the basic functions and operators used in spatial statistics
may be represented using rank-structured tensor formats and that the error of this 
representation exhibits the exponential decay with respect to the tensor rank.
We applied the Tucker and canonical tensor decompositions to a family of 
Mat\'ern- and Slater-type functions with varying parameters and demonstrated numerically 
that their approximations exhibit exponentially fast convergence. 
A low-rank tensor approximation of the Mat\'{e}rn covariance function and its Fourier 
transform is considered. We separated the radial basis functions using the Laplace transforms 
to prove the existence of such low-rank approximations, and applied the $\sinc$ quadrature method 
to estimate the tensor ranks and accuracy.

We also demonstrated how to compute $\diag(\bC)$, $\trace(\bC)$, the matrix-vector product, 
Kriging operations, and the geostatistical optimal design in a low-rank tensor format with at a linear cost.
 {
For matrix $\bC$ of size $N\times N$, $N=6000^3$ and of tensor rank 1, 
we were able to compute the Cholesky factorization in 15 sec. We also computed the 
Tucker approximation to the 3D Slater function $e^{-\Vert x\Vert}$ on the grid 
with $513^{3}$ points and Tucker ranks $r=1,2,...,10$. Furthermore, 
we demonstrated how to compute $\trace(\bC)$ for $N=n^d=1000^{1000}$. 
This might be useful in machine learning.}
  In this paper,  
operations such as computing the Cholesky factorization, inverse, and determinant 
  have been implemented   
for rank-1 tensors (e.g., the Gaussian covariance has a tensor rank-1). 
These formulas  could be useful for developing successive 
rank-1 updates in greedy algorithms. 
Further investigations are needed for the representation of these quantities with the ranks higher than one. 

Additionally, in Section~\ref{ssec:Numerics1} we studies the influence of the 
parameters of the Mat\'{e}rn covariance function on the tensor ranks 
(Figures \ref{fig:Conv_Tucker_allp} and \ref{fig:Conv_Tucker_mg}). 
We observed (see Fig.~\ref{fig:Tuck_SD_Matern}) that the dependence 
of the parameters of the Mat\'{e}rn covariance function on the tensor ranks
is very weak, and the ranks grew slowly. 
In this paper, we also highlighted that big data statistical problems can be effectively treated by 
using the special low-rank tensor techniques.

\section*{Acknowledgements}
The research reported in this publication was supported by funding from King 
Abdullah University of Science and Technology (KAUST).

  \bibliographystyle{siam}
\begin{footnotesize}
\bibliography{Hcovariance2}
 \end{footnotesize}



%
\end{document}